\documentclass[a4paper,11pt]{amsart}
\usepackage[utf8x]{inputenc}
\usepackage{amsthm,amsfonts,amsmath,amssymb,enumerate}

\usepackage{pstricks-add}

\usepackage{pgf,tikz}
\usetikzlibrary{arrows}

\usepackage{pdfsync}
\newcommand{\N}{{\mathbb N}}

\newcommand{\K}{{\mathbb K}}
\renewcommand{\a}{{\mathfrak a}}
\newcommand{\n}{{\mathfrak n}}

\newcommand{\g}{{\mathfrak g}}
\newcommand{\s}{{\mathfrak s}}
\newcommand{\h}{{\mathfrak h}}
\newcommand{\f}{{\mathfrak{f}}}
\renewcommand{\b}{{\mathfrak{b}}}
\newcommand{\gl}{{\mathfrak{gl}}}
\renewcommand{\sl}{{\mathfrak{sl}}}

\newcommand{\End}{\operatorname{End}}

\newcommand{\Der}{\operatorname{Der}}
\newcommand{\rank}{\operatorname{rank}}
\newcommand{\HomL}{\operatorname{Hom_{\text{Lie}}}}
\newcommand{\ad}{\operatorname{ad}}

\renewcommand{\Im}{\operatorname{Im}}
\newcommand{\Hom}{\operatorname{Hom}}
\newcommand{\Aut}{\operatorname{Aut}}
\newcommand{\Id}{\operatorname{Id}}
\newcommand{\nilsh}{\operatorname{nilsh}}

\newtheorem{theorem}{Theorem}[section]
\newtheorem{lemma}[theorem]{Lemma}
\newtheorem{corollary}[theorem]{Corollary}

\newtheorem{proposition}[theorem]{Proposition}

\theoremstyle{definition}
\newtheorem{definition}[theorem]{Definition}

\theoremstyle{remark}
\newtheorem{remark}[theorem]{Remark}

\numberwithin{equation}{section}
\begin{document}

\title[Cohomology of solvable Lie algebras and deformations]{Total cohomology of solvable Lie algebras \\ and linear deformations}
\author{Leandro Cagliero and Paulo Tirao}
\address{CIEM-FaMAF, Universidad Nacional de C\'ordoba, Argentina}
\date{November 29, 2012}

\maketitle

\begin{abstract}
Given a finite dimensional Lie algebra $\g$, let 
$\Gamma_\circ(\g)$ be the set of irreducible $\g$-modules 
with non-vanishing cohomology. 
We prove that a $\g$-module $V$ belongs to $\Gamma_\circ(\g)$ only if $V$ is contained  
in the exterior algebra of the solvable radical $\s$ of $\g$, 
showing in particular that $\Gamma_\circ(\g)$ is a finite set and we deduce
that $H^*(\g,V)$ is an $L$-module, where $L$ is a fixed subgroup of 
the connected component of $\Aut(\g)$ which contains a Levi factor.

We describe $\Gamma_\circ$ in some basic examples, including the Borel subalgebras,
and we also determine $\Gamma_\circ(\s_n)$ for an extension $\s_n$ of the 2-dimensional abelian Lie algebra 
by the standard filiform Lie algebra $\f_n$.
To this end, we described the cohomology of $\f_n$.

We introduce the \emph{total cohomology} of a Lie algebra $\g$, as 
$TH^*(\g)=\bigoplus_{V\in \Gamma_\circ(\g)} H^*(\g,V)$ and
we develop further the theory of linear deformations in order to prove that the total 
cohomology of a solvable Lie algebra is the cohomology of its nilpotent shadow.
Actually we prove that $\s$ lies, in the variety of Lie algebras, 
in a linear subspace of dimension at least $\dim (\s/\n)^2$, $\n$ being the nilradical of $\s$, 
that contains the nilshadow of $\s$ and such that all its points have the same total cohomology.
\end{abstract}

\section{Introduction}

In this paper $\K$ will be a field of characteristic 0.
In addition, since many of our results are simpler to present having available
eigenvalues, we will assume that $\K$ is algebraically closed. 
Nevertheless, as it is usual when dealing with cohomology, 
many of our results are true (perhaps with minor modifications) 
without this assumption.
All Lie algebras and modules will be of finite dimension over $\K$.

A well known result of Whitehead (see for instance \cite{J}) 
states that the cohomology of a semisimple Lie algebra $\g$ 
with coefficients in a $\g$-module $V$ vanishes if and only if 
$V$ does not contain the trivial module.
The same result for nilpotent Lie algebras was proved by Dixmier \cite{D}.
It is natural to consider the same problem for more general families of Lie algebras $\g$, 
that is, determining the set $\Gamma(\g)$ of modules with non-vanishing cohomology.

Very little is known about this problem.
A natural first step is to determine the set
$\Gamma_\circ(\g)$ of irreducible modules with non-vanishing cohomology,
but also in this case almost nothing is known.
Already in the solvable case
the determination of $\Gamma_\circ$, that is the set of characters 
with non-vanishing cohomology, is an interesting problem that seems to be quite difficult.

The Hochschild-Serre spectral sequence in \cite{HS} describes the cohomology of a Lie algebra $\g$
in terms of the cohomologies of an ideal $\h$ and that of the quotient Lie algebra $\g/\h$.
This description is in general sophisticated,
however it is extremely clear in some particular instances.
For example, in the case when $\h=\s$ is the solvable radical of $\g$, one has the remarkable formula
\begin{equation}\label{eqn:hochschild-serre}
 H^*(\g,V) \simeq H^*(\g/\s,\K)\otimes H^*(\s,V)^\g.
\end{equation}
Whitehead's result describing $\Gamma(\g)$ for a semisimple $\g$ follows at once from this formula. 
The main tool for Dixmier's result describing $\Gamma(\n)$ for a nilpotent $\n$ also
follows from the Hochschild-Serre spectral sequence.

In this paper we investigate $\Gamma_\circ(\g)$. 
The results of Hochschild-Serre play a key role in different ways.
Formula \eqref{eqn:hochschild-serre}
shows that the solvable case is a significant obstruction to understand the general case
and we focus mainly on this case.
Some of our results are presented in the
framework of deformation theory of Lie algebras.

\medskip

In the first part of the paper we prove that, 
for an irreducible $\g$-module $V$, $H^*(\g,V)$ vanishes unless $V$ is contained 
(is a submodule of a quotient module) 
in the exterior algebra of the solvable radical $\s$ of $\g$, 
showing in particular that $\Gamma_\circ(\g)$ is a finite set. 

In \S \ref{sec:examples} we describe  $\Gamma_\circ$ in 
some basic examples, including the Borel subalgebras, that illustrate 
some aspects of the problem of determining $\Gamma_\circ$.
We also determine $\Gamma_\circ(\s_n)$ in a more involved case, 
that of an extension $\s_n$ of the 2-dimensional abelian Lie algebra 
by the standard filiform Lie algebra $\f_n$.
To this end, we described the cohomology of $\f_n$.
We close this section showing that, for a solvable Lie algebra $\s$, 
$\Gamma_\circ(\s)$ is contained in the set of $\s$-characters 
of $\s/\s''$.

In \S \ref{sec:total} we introduce the \emph{total cohomology} $TH^*(\g)$ of a Lie algebra $\g$, as 
\[ TH^*(\g)=\bigoplus_{V\in \Gamma_\circ(\g)} H^*(\g,V). \]
We show that $TH^*(\g)$ is an $L$-module where $L$ is a certain subgroup of 
the connected component of $\Aut(\g)$ which contains a Levi factor. 
Actually $H^*(\g,V)$ is an $L$-submodule of $TH^*(\g)$ for all $V$.

\medskip
 
In the second part of the paper we investigate the total cohomology of solvable Lie algebras.
We start by proving that if $\s$ is a solvable Lie algebra, 
then there is another solvable algebra $\tilde \s$, that is `more nilpotent' than $\s$ (the nilradical grows),
such that $TH^*(\s)=TH^*(\tilde\s)$. 
An induction argument leads to a nilpotent Lie algebra $\n$ whose cohomology
is isomorphic to $TH^*(\s)$.
This process can be understood in the framework of (linear) deformations of Lie algebras.
We develop further the theory of linear deformations in order to prove that the total 
cohomology of a solvable Lie algebra is the
cohomology of its nilpotent shadow or simply nilshadow.
Actually we prove that $\s$ lies, in the variety of Lie algebras, 
in a large linear subspace (of dimension at least $\dim (\s/\n)^2$, $\n$ being the nilradical of $\s$) 
that contains the nilshadow of $\s$ and such that all its points have the same total cohomology.

\section{Preliminaries}\label{sec:preliminaries}
In this section we describe explicitly some basic results 
obtained from the Hochschild-Serre spectral sequence in 
two particular instances: when the ideal $\h$ has codimension 1, 
and when $\g=\a \ltimes \h$, where
$\a$ is an abelian subalgebra of semisimple derivations of $\h$.

Let $\g$ be a Lie algebra, $\h$ an ideal of codimension one, and let $D$ be an element outside $\h$.
Let $V$ be a $\g$-module and $\rho$ be the corresponding representation.
Then there is an exact sequence
\begin{equation}\label{eqn:long-exact-seq}
\dots\rightarrow H^p(\g,V) \rightarrow H^p(\h,V)
\overset{\theta_\rho(D)}{\longrightarrow} H^p(\h,V) \rightarrow H^{p+1}(\g,V)\rightarrow \dots
\end{equation}
where the morphism $\theta_\rho(D)$ is induced by the \emph{Lie derivative} of $D$.
Recall that, if $f\in \Lambda \h^* \otimes V$, then
\[
\theta_\rho(D)(f)(x)= -f(\overline{D}(x)) + \rho(D)f(x), \quad x\in\Lambda \h,
\]
where $\overline{D}$ is the canonical extension of $\ad D$ to the exterior algebra of $\h$.
Denoting by
\[ \theta_0(D)=-\overline{D}^t \otimes \Id \]
we have $\theta_\rho(D) = \theta_0(D) + \Id\otimes \rho(D)$ and for simplicity we write
\[ \theta_\rho(D) = \theta_0(D) + \rho(D). \]
Since (\ref{eqn:long-exact-seq}) is exact, it follows that, for every $p=0,\dots, \dim \g$,  
\begin{align}\label{eqn:dixmier1}
 H^p(\g,V) &= H^p(\h,V)^{\theta_\rho(D)} \oplus H^{p-1}(\h,V)_{\theta_\rho(D)}\otimes D^*, \\ \label{eqn:dixmier2}
           &\simeq H^p(\h,V)^{\theta_\rho(D)} \oplus H^{p-1}(\h,V)_{\theta_\rho(D)},
\end{align}
where $D^*\in\g^*$ is the linear functional that vanishes on $\h$ dual to $D$.

The space of coinvariants $H^{p-1}(\h,V)_{\theta_\rho(D)}$ has 
the same dimension as the space of invariants $H^{p-1}(\h,V)^{\theta_\rho(D)}$.
Moreover, if $\theta_\rho(D)$ is diagonalizable (which is the case if 
$\ad D$ and $\rho(D)$ are diagonalizable) then 
we can identify canonically the invariants and coinvariants and thus 
\begin{equation}\label{eqn:dixmier3}
 H^p(\g,V) = H^p(\h,V)^{\theta_\rho(D)} \oplus H^{p-1}(\h,V)^{\theta_\rho(D)}\otimes D^*.
\end{equation}
An induction argument proves the following Hochschild-Serre type formula.
\begin{proposition}\label{prop:hoch-serre-abelian}
Let $\g=\a \ltimes \h$, where $\a$ is an abelian subalgebra of semisimple derivations 
of the ideal $\h$.
Let $V$ be a $\g$-module such that $\rho(A)$ is diagonalizable for all $A\in\a$.
Then 
\[
   H^*(\g,V) \simeq H^*(\h,V)^\a\otimes \Lambda \a^*,
\]
where $\a$ acts on $H^*(\h,V)$ by $\theta_\rho$.
In particular, if $V$ is a 1-dimensional module corresponding 
to a character $\lambda$ and $\h\subset\ker\lambda$, then 
\[
   H^*(\g,\lambda) \simeq  H^*(\h)_{(-\lambda)}\otimes\Lambda \a^*,
\]
where $H^*(\h)_{(-\lambda)}$ is the $\a$-isotopic component of type $-\lambda$
where $\a$ acts on $H^*(\h,V)$ by $\theta_0$.
\end{proposition}

\section{Modules with non-vanishing cohomology}

Given a Lie algebra $\g$ over $\mathbb K$, 
denote by $\text{Rep}(\g)$ and $\text{Irrep}(\g)$ respectively the sets of equivalence classes of finite dimensional
representations and finite dimensional irreducible representations of $\g$.
If $(\tau,W)\in\text{Rep}(\g)$, denote by  
\[
 \chi(\g,W)=\left\{(\pi,V)\in\text{Irrep}(\g)
\left|
\begin{matrix}
 (\pi,V)\text{ is an irreducible quotient}\\
\text{ in a composition series of $W$}.
\end{matrix}
\right.
\right\}
\]
By the Jordan-H\"older theorem,  $\chi(\g,W)$ is well defined.  

If $\g=\s$ is solvable the irreducible representations of $\s$
are all one dimensional and $\chi(\s,W) \subseteq(\s/\s')^*$ for any representation $(\pi,W)$.
Moreover, from Lie's theorem it follows that $\chi(\s,W)$ is the set of diagonal matrix entries 
in a triangularized form of $\pi$.
For example
\begin{equation}\label{eqn:chi-solvable}
 \chi(\s,\Lambda^p\s)=\left\{\sum_{i=1}^p\lambda_i:\lambda_i\in\chi(\s,\s)\right\}
\end{equation}
and, in particular, $\chi(\n,\Lambda^p\n)=0$ if $\n$ is nilpotent.

\

If  $(\tau_1,W_1)$ and $(\tau_2,W_2)$ are two representations of a Lie algebra $\g$, then 
the basic results about the Grothendieck ring of $\g$ (see for instance \cite{CR})
imply that 
\[
\chi(\g,W_1 \otimes W_2)=\bigcup_{V_1\in \chi(\g,W_1)}
\bigcup_{V_2\in \chi(\g,W_2)}\chi(\g,V_1\otimes V_2).
\]
In turn, it is a classical result of Chevalley that for any $\g$,
the tensor product of irreducible representations of $\g$ is completely reducible.
In particular, if  $(\pi_i,V_i)\in \text{Irrep}(\g)$, $i=1,2,3$, then 
$(\pi_1^*,V_1^*)\in\chi(\g,V_2\otimes V_3)$ 
if and only if $(V_1\otimes V_2\otimes V_3)^{\g}\ne0$.
This observation proves the following reciprocity result.

\begin{proposition}\label{prop.reciprocity}
Let $(\pi_1,V_1)$, $(\pi_2,V_2)$ and $(\tau,W)$ be representations of $\g$, 
$(\pi_1,V_1)$ and $(\pi_2,V_2)$ being irreducible.
Then  
\[
(\pi_1,V_1)\in \chi(\g,W^*\otimes V_2^*) \text{ if and only if } (\pi_2,V_2)\in \chi(\g,W^*\otimes V_1^*). 
\]
In particular, the trivial representation of $\g$ belongs to $\chi(\g,W^*\otimes V_1)$
if and only if  $(\pi_1,V_1)$ belongs to $\chi(\g,W)$.
\end{proposition}

\

Let us now consider the family of $\g$-modules with non-vanishing cohomology  
distinguishing the irreducible ones.
For $p=0,\dots,\dim\g$, let
\begin{eqnarray*}
 \Gamma^p(\g) &=& \{(\pi,W)\in\text{Rep}(\g):H^p(\g,W)\ne0\}, \\
 \Gamma_\circ^p(\g) &=& \{(\pi,V)\in\text{Irrep}(\g):H^p(\g,V)\ne0\}.
\end{eqnarray*}
Also denote
\[ \Gamma(\g)=\bigcup_p \Gamma^p(\g) \qquad\text{and}\qquad \Gamma_\circ(\g)=\bigcup_p \Gamma^p_\circ(\g). \]

It is well known that $\Gamma_\circ^p(\g)=0$ for all $p$ if $\g$ 
is semisimple \cite[Ch. III, Thm.\ 14]{J} or nilpotent \cite{D} and
moreover, in these two cases, $(\pi,W)\in\Gamma^p(\g)$ if and only if the trivial representation is contained in $W$. 
We now prove the following more general result.

\begin{theorem}\label{thm:TH} 
Let $\g=\g_0\ltimes\s$ be the Levi decomposition of a 
finite dimensional Lie algebra $\g$. 
Then
\[
 \Gamma_\circ^p(\g)\subseteq\chi(\g,\Lambda^p \s),
\]
for all $0\le p\le\dim\g$. 
In particular $\Gamma_\circ^p(\g)$ is a finite set for all $p$.
\end{theorem}

\begin{proof}
From the Hochschild-Serre formula \eqref{eqn:hochschild-serre}
\[
 H^p(\g,V)=\sum_{p_1+p_2=p}H^{p_1}(\g_0)\otimes H^{p_2}(\s,V)^{\g}
\]
it follows that if $(\pi,V)\in\Gamma_\circ^p(\g)$, then 
there exists a non-zero subquotient of $(\Lambda^p\s)^*\otimes V$
that is trivial as a $\g$-module.
This implies that the trivial representation of $\g$ belongs to 
$\chi(\g,(\Lambda^p\s)^*\otimes V)$ and, by Proposition \ref{prop.reciprocity},
we conclude that $(\pi,V)\in\chi(\g,\Lambda^p\s)$.
\end{proof}

\begin{remark}
The above theorem gives a first estimate for $\Gamma^p_\circ(\g)$.
Recall that Whitehead's and Dixmier's theorems state that a $\g$-module $V$ is in $\Gamma(\g)$
if and only if $V$ contains the (only) irreducible $\g$-module in $\Gamma_\circ(\g)=\{ \text{trivial module}\}$
for $\g$ semisimple or nilpotent.
This is far from being true in general.
It seems to be much more difficult to obtain an estimation for $\Gamma^p(\g)$.
For instance, if $\b$ is a Borel subalgebra of a semisimple Lie algebra, 
then the adjoint module $\b$ does not
belong to $\Gamma(\b)$ \cite{LL},
however $\b$ contains many irreducible representations in $\Gamma_\circ(\b)$
as we show in the next section. 
\end{remark}

\section{The solvable case: characters with non-vanishing cohomology}\label{sec:examples}

The problems of understanding the sets $\Gamma^p(\s)$ and $\Gamma_\circ^p(\s)$ are difficult
already for solvable Lie algebras $\s$.
If $\s$ is solvable, it follows from Theorem \ref{thm:TH} and (\ref{eqn:chi-solvable}) that
\[  
 \Gamma_\circ^p(\s)\subseteq\chi(\s,\Lambda^p \s)=\left\{\sum_{i=1}^p\lambda_i:\lambda_i\in\chi(\s,\s)\right\}.
\]
The inclusion is in general proper, but not always.
In the following examples we will compute $\Gamma_\circ^p(\s)$ in some particular cases
that exhibit a variety of situations that may occur. 
At the end of this section we prove that
\[
\Gamma_\circ^1(\s)\subset\chi(\s, \s/\s''),
\]
where $\s''$ is the second term of the derived series of $\s$.

\subsection{A first elementary example}

This is an elementary example of a solvable Lie algebra $\s$ and a character $\lambda\in\chi(\s,\Lambda^p\s)$
such that $H^p(\s,\lambda)=0$, for all $p$.

Let $\s$ be the 4-dimensional solvable Lie algebra with basis $\{ u,x,y,z\}$ and bracket defined by:
 \[ [u,x]=x, [u,y]=y, [u,z]=2z, [x,y]=z. \]
All characters of $\s$ are $\lambda_a$ for some $a\in\K$, where 
 \[ \lambda_a(u)=a, \quad \lambda_a(x)=\lambda_a(y)=\lambda_a(z)=0. \]
The characters in the adjoint representation of $\s$ are clearly $\lambda_0$, $\lambda_1$ and $\lambda_2$.
Therefore
\begin{eqnarray*}
 \chi(\s,\Lambda^0\s) &=& \{\lambda_0\} \\
 \chi(\s,\Lambda^1\s) &=& \{\lambda_0,\lambda_1,\lambda_2\} \\
 \chi(\s,\Lambda^2\s) &=& \{\lambda_1,\lambda_2,\lambda_3\} \\
 \chi(\s,\Lambda^3\s) &=& \{\lambda_2,\lambda_3,\lambda_4\} \\
 \chi(\s,\Lambda^4\s) &=& \{\lambda_4\}. \\
\end{eqnarray*}
From Theorem \ref{thm:TH} it follows that $H^*(\s,\lambda_a)$ vanishes if $a\not\in\{0,1,2,3,4\}$.
The following table, obtained by direct calculations, shows the cohomology $H^p(\s,\lambda_a)$ of $\s$ with coefficients
in the character modules $\lambda_a$ with $a=0,1,2,3,4$ and for all $p=0,\dots,4$.
The entries `$-$' indicate that, by Theorem \ref{thm:TH}, there is no possible cohomology there.

\begin{center}
 \begin{tabular}{c|ccccc}
       & $a=0$ & $a=1$ & $a=2$ & $a=3$ & $a=4$ \\ \hline
 $p=0$ & 1 & $-$ & $-$ & $-$ & $-$ \\
 $p=1$ & 1 &  2 & 0 & $-$ & $-$ \\
 $p=2$ & $-$ & 2 & 0 & 2 & $-$ \\
 $p=3$ & $-$ & $-$ & 0 & 2 & 1 \\
 $p=4$ & $-$ & $-$ & $-$ & $-$ & 1 \\
 \end{tabular}

\end{center}
Therefore $\Gamma_\circ^0(\s)=\{ \lambda_0\}$, $\Gamma_\circ^1(\s)=\{ \lambda_0,\lambda_1\}$, 
$\Gamma_\circ^2(\s)=\{ \lambda_1,\lambda_3\}$, $\Gamma_\circ^3(\s)=\{ \lambda_3,\lambda_4\}$ and 
$\Gamma_\circ^4(\s)=\{\lambda_4\}$.

\subsection{A second elementary example}

Here we show an example where $\Gamma_\circ=\chi(\s,\Lambda\s)$.
Let $\h$ be an abelian Lie algebra and $D$ any linear operator on $\h$.
Consider the solvable Lie algebra $\s=D\ltimes\h$.
Given a character $\lambda$ of $\s$, it follows from \eqref{eqn:dixmier2} that
\[
\dim H^p(\s,\lambda)=\dim \big(\Lambda^p\h^*\big)^{\theta(D)+\lambda(D)}+\dim \big(\Lambda^{p-1}\h^*\big)^{\theta(D)+\lambda(D)}.
\]
Identifying $\chi(\s,\Lambda^p\s)$ with the set of eigenvalues of $\overline{D}$ on $\Lambda^p\h$,
it follows that $\lambda\in\chi(\s,\Lambda^p\h)$ or $\lambda\in\chi(\s,\Lambda^{p-1}\h)$ if and only if  $H^p(\s,\lambda)\ne 0$.
Therefore $\Gamma_\circ^p(\s)=\chi(\s,\Lambda^p\h)\cup \chi(\s,\Lambda^{p-1}\h)$ and in particular
$\Gamma_\circ(\s)=\chi(\s,\Lambda\s)$.

\subsection{A third elementary example}

Let $\h$ be an abelian Lie algebra of dimension $n$, let $\a$ be an abelian subalgebra of diagonalizable operators on $\h$ of dimension $r$
and consider the solvable Lie algebra $\s= \a \ltimes \h$.

On the one hand the set $\chi(\s,\Lambda^p\s)$ can be identified with the set of weights of $\a$ on $\Lambda^p\h$.
On the other hand, by Proposition \ref{prop:hoch-serre-abelian}
\[
   H^*(\s,\lambda) \simeq \big(\Lambda\h^*\big)_{(-\lambda)} \otimes \Lambda \a^*.
\]
It follows that for $p=0,\dots,n+r$
\[ 
\Gamma_\circ^p(\s) = \chi(\s,\Lambda^p\h)\cup \chi(\s,\Lambda^{p-1}\h) \cup \dots \cup \chi(\s,\Lambda^{p-r}\h).
 \]
Again, as in the previous example, $\Gamma_\circ(\s)=\chi(\s,\Lambda\s)$.

\subsection{Borel subalgebras}

Let $\b=\h+\n$ be a Borel subalgebra of a semisimple Lie algebra $\g$
of rank $r$, and let $\rho$ be the semisum of all positive roots.
Let $W$ be the Weyl group of $\g$ and, for $w\in W$, let $\ell(w)$ be the length of $w$.
 
If $\lambda$ is a weight of $\b$, 
it follows from Proposition \ref{prop:hoch-serre-abelian} that
\[
H^p(\b,\lambda)=\sum H^{p-i}(\n)_{(-\lambda)}\otimes \Lambda^i \h^*,
\]
where $H^j(\n)_{(-\lambda)}$ is the $-\lambda$ component of $H^j(\n)$ as $\h$-module. 
It follows from Kostant's description of the cohomology of nilradicals of parabolic
subalgebras that
\[ 
\Gamma_\circ^p(\b)=\bigcup_{j=p-r}^p 
\{\lambda\in\chi(\b,\Lambda\n):\lambda|_{\h}=\rho-w\rho
\text{ for some $w\in W$ with $\ell(w)=j$}\}. 
\]
We know that, after restricting to $\h$, $\chi(\b,\Lambda\b)$ is the convex (root-lattice) 
polytope whose points are all the sums of distinct roots. 
Therefore we have seen that, in this case, $\Gamma_\circ^p$ is precisely the extreme vertices of 
the polytope $\chi(\b,\Lambda\b)$.

The sets $\chi(\b,\Lambda \b)$ and $\chi(\b,\Lambda \n)$ coincide
(their multiplicities are related by a factor $2^r$).
The following diagrams show the sets of weights $\chi(\b,\Lambda\n)$ with multiplicities and marked with a bullet
those in $\Gamma_\circ(\b)$ for the Borel subalgebras $\b$ of $\g$ of type $A_2$ and $G_2$.
The arrows indicated the positive roots of $\g$.

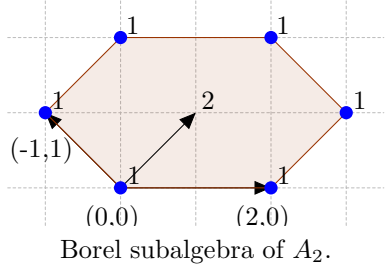
\begin{center}


\small
\definecolor{qqqqqq}{rgb}{0,0,0}
\definecolor{zzttqq}{rgb}{0.6,0.2,0}
\definecolor{qqqqff}{rgb}{0,0,1}
\definecolor{cqcqcq}{rgb}{0.75,0.75,0.75}
\begin{tikzpicture}[line cap=round,line join=round,>=triangle 45,x=1.0cm,y=1.0cm]
\draw [color=cqcqcq,dash pattern=on 1pt off 1pt, xstep=1.0cm,ystep=1.0cm] (-1.5,-0.5) grid (3.5,2.5);
\clip(-1.5,-0.5) rectangle (3.5,2.5);
\fill[color=zzttqq,fill=zzttqq,fill opacity=0.1] (0,0) -- (-1,1) -- (0,2) -- (2,2) -- (3,1) -- (2,0) -- cycle;
\draw [->,color=qqqqqq] (0,0) -- (2,0);
\draw [->,color=qqqqqq] (0,0) -- (-1,1);
\draw [->,color=qqqqqq] (0,0) -- (1,1);
\draw [color=zzttqq] (0,0)-- (-1,1);
\draw [color=zzttqq] (-1,1)-- (0,2);
\draw [color=zzttqq] (0,2)-- (2,2);
\draw [color=zzttqq] (2,2)-- (3,1);
\draw [color=zzttqq] (3,1)-- (2,0);
\draw [color=zzttqq] (2,0)-- (0,0);
\draw [color=qqqqqq](-0.05,2.4) node[anchor=north west] {1};
\draw [color=qqqqqq](1.95,2.4) node[anchor=north west] {1};
\draw [color=qqqqqq](2.97,1.4) node[anchor=north west] {1};
\draw [color=qqqqqq](1.96,0.4) node[anchor=north west] {1};
\draw [color=qqqqqq](-0.03,0.4) node[anchor=north west] {1};
\draw [color=qqqqqq](-1.04,1.4) node[anchor=north west] {1};
\draw [color=qqqqqq](0.95,1.4) node[anchor=north west] {2};
\draw [color=qqqqqq](-0.6,-0.1) node[anchor=north west] {(0,0)};
\draw [color=qqqqqq](1.4,-0.1) node[anchor=north west] {(2,0)};
\draw [color=qqqqqq](-1.6,0.8) node[anchor=north west] {(-1,1)};
\fill [color=qqqqff] (0,0) circle (2.5pt);
\fill [color=qqqqff] (2,0) circle (2.5pt);
\fill [color=qqqqff] (-1,1) circle (2.5pt);
\fill [color=qqqqff] (0,2) circle (2.5pt);
\fill [color=qqqqff] (2,2) circle (2.5pt);
\fill [color=qqqqff] (3,1) circle (2.5pt);
\end{tikzpicture}
 \\
Borel subalgebra of $A_2$.

\ \\

\ \\


\small
\definecolor{qqqqqq}{rgb}{0,0,0}
\definecolor{zzttqq}{rgb}{0.6,0.2,0}
\definecolor{qqqqff}{rgb}{0,0,1}
\definecolor{cqcqcq}{rgb}{0.75,0.75,0.75}
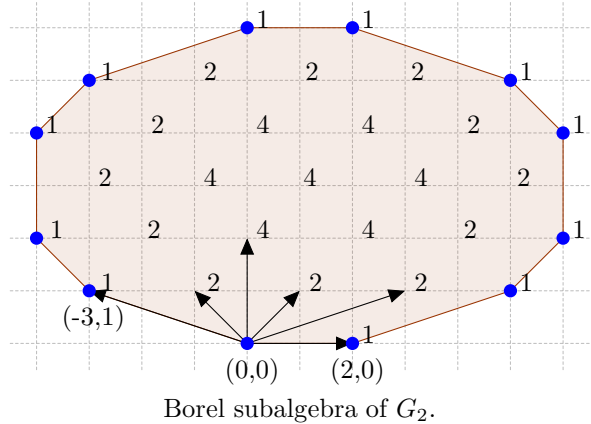
\begin{tikzpicture}[line cap=round,line join=round,>=triangle 45,x=0.7cm,y=0.7cm]
\draw [color=cqcqcq,dash pattern=on 1pt off 1pt, xstep=0.7cm,ystep=0.7cm] (-4.5,-0.5) grid (6.5,6.5);
\clip(-4.5,-0.8) rectangle (6.5,6.5);
\fill[color=zzttqq,fill=zzttqq,fill opacity=0.1] (0,0) -- (-3,1) -- (-4,2) -- (-4,4) -- (-3,5) -- (0,6) -- (2,6) -- (5,5) -- (6,4) -- (6,2) -- (5,1) -- (2,0) -- cycle;
\draw [color=zzttqq] (0,0)-- (-3,1);
\draw [color=zzttqq] (-3,1)-- (-4,2);
\draw [color=zzttqq] (-4,2)-- (-4,4);
\draw [color=zzttqq] (-4,4)-- (-3,5);
\draw [color=zzttqq] (-3,5)-- (0,6);
\draw [color=zzttqq] (0,6)-- (2,6);
\draw [color=zzttqq] (2,6)-- (5,5);
\draw [color=zzttqq] (5,5)-- (6,4);
\draw [color=zzttqq] (6,4)-- (6,2);
\draw [color=zzttqq] (6,2)-- (5,1);
\draw [color=zzttqq] (5,1)-- (2,0);
\draw [color=zzttqq] (2,0)-- (0,0);
\draw [color=qqqqqq](0.0,6.5) node[anchor=north west] {1};
\draw [color=qqqqqq](2.0,6.5) node[anchor=north west] {1};
\draw [color=qqqqqq](5.0,5.5) node[anchor=north west] {1};
\draw [color=qqqqqq](6.0,4.5) node[anchor=north west] {1};
\draw [color=qqqqqq](6.0,2.5) node[anchor=north west] {1};
\draw [color=qqqqqq](5.0,1.5) node[anchor=north west] {1};
\draw [color=qqqqqq](2.0,0.5) node[anchor=north west] {1};
\draw [color=qqqqqq](-2.95,1.5) node[anchor=north west] {1};
\draw [color=qqqqqq](-3.92,2.5) node[anchor=north west] {1};
\draw [color=qqqqqq](-4,4.5) node[anchor=north west] {1};
\draw [color=qqqqqq](-2.95,5.5) node[anchor=north west] {1};
\draw [color=qqqqqq](0.0,2.5) node[anchor=north west] {4};
\draw [color=qqqqqq](-2.07,2.5) node[anchor=north west] {2};
\draw [color=qqqqqq](3.0,1.5) node[anchor=north west] {2};
\draw [color=qqqqqq](-0.94,1.5) node[anchor=north west] {2};
\draw [color=qqqqqq](1.0,1.5) node[anchor=north west] {2};
\draw [color=qqqqqq](3.94,2.5) node[anchor=north west] {2};
\draw [color=qqqqqq](2.0,2.5) node[anchor=north west] {4};
\draw [color=qqqqqq](-1.0,5.5) node[anchor=north west] {2};
\draw [color=qqqqqq](2.94,5.5) node[anchor=north west] {2};
\draw [color=qqqqqq](-0.0,4.5) node[anchor=north west] {4};
\draw [color=qqqqqq](-2.0,4.5) node[anchor=north west] {2};
\draw [color=qqqqqq](4.0,4.5) node[anchor=north west] {2};
\draw [color=qqqqqq](2.0,4.5) node[anchor=north west] {4};
\draw [color=qqqqqq](0.89,3.5) node[anchor=north west] {4};
\draw [color=qqqqqq](-3.0,3.5) node[anchor=north west] {2};
\draw [color=qqqqqq](4.95,3.5) node[anchor=north west] {2};
\draw [color=qqqqqq](3,3.5) node[anchor=north west] {4};
\draw [color=qqqqqq](-1.0,3.5) node[anchor=north west] {4};
\draw [color=qqqqqq](0.93,5.5) node[anchor=north west] {2};
\draw [color=qqqqqq](-0.6,-0.1) node[anchor=north west] {(0,0)};
\draw [color=qqqqqq](1.4,-0.1) node[anchor=north west] {(2,0)};
\draw [color=qqqqqq](-3.7,0.9) node[anchor=north west] {(-3,1)};
\draw [->,color=qqqqqq] (0,0) -- (2,0);
\draw [->,color=qqqqqq] (0,0) -- (-1,1);
\draw [->,color=qqqqqq] (0,0) -- (1,1);
\draw [->,color=qqqqqq] (0,0) -- (0,2);
\draw [->,color=qqqqqq] (0,0) -- (-3,1);
\draw [->,color=qqqqqq] (0,0) -- (3,1);
\begin{scriptsize}
\fill [color=qqqqff] (0,0) circle (2.5pt);
\fill [color=qqqqff] (2,0) circle (2.5pt);
\fill [color=qqqqff] (-3,1) circle (2.5pt);
\fill [color=qqqqff] (-4,2) circle (2.5pt);
\fill [color=qqqqff] (-4,4) circle (2.5pt);
\fill [color=qqqqff] (-3,5) circle (2.5pt);
\fill [color=qqqqff] (0,6) circle (2.5pt);
\fill [color=qqqqff] (2,6) circle (2.5pt);
\fill [color=qqqqff] (5,5) circle (2.5pt);
\fill [color=qqqqff] (6,4) circle (2.5pt);
\fill [color=qqqqff] (6,2) circle (2.5pt);
\fill [color=qqqqff] (5,1) circle (2.5pt);
\end{scriptsize}
\end{tikzpicture}
 \\
Borel subalgebra of $G_2$.
\end{center}

\subsection{The standard filiform Lie algebra $\f_n$ and a solvable Lie algebra with nilradical $\f_n$}

Let $\mathfrak{f}_n$ be the standard filiform Lie algebra of dimension $n+2$. 
It is not difficult to see that $\text{Der}(\f_n)$ has a maximal toral subalgebra $\a$ 
of dimension 2 (its elements are ad-semisimple).
In this example we describe the $\a$-module structure of the cohomology of $\mathfrak{f}_n$ 
and we obtain the sets $\Gamma_\circ^p(\s_n)$ for the solvable Lie algebra $\s_n=\a\ltimes\f_n$.

Given $n\in\N$, let $V_n$ be the irreducible representation of $\sl(2)$ 
with highest weight $n$ and let $\{e_{n-2k}:k=0,\dots,n\}$ be a basis of weight vectors of $V_n$. 
Let $\{E,H,F\}$ be the standard basis of $\sl(2)$ and let $I$ be the  identity of $\gl(2)$.  
Consider $V_n$ as a representation of $\gl(2)$ by letting $I$ act as the identity map on $V_n$.
We have  
 \begin{align*}
  \mathfrak{f}_n &\simeq E\ltimes V_n, \\
   \a&\simeq\text{span}\{H,I\}, \\
  \s_n = \a\ltimes\mathfrak{f}_n &\simeq \text{span}\{H,I,E\}\ltimes V_n,
 \end{align*}
where the right hand sides are considered inside $\gl(2)\ltimes V_n$.
 It is clear that 
$\a$ is a Cartan subalgebra of $\s_n$ and $\s_n'= \mathfrak{f}_n$.
Therefore we can identify any weight $\lambda$ of $\s_n$ with the pair $(\lambda(H),\lambda(I))$.
In particular, the set of roots of $\s_n$ with respect to $\a$ is 
\[
 \Phi(\s_n,\a)=\{\alpha=(2,0),\beta_k=(n-2k,1):\;k=0,\dots,n\}
\]
and $\chi(\s_n,\Lambda \s_n)$ is set of root-lattice points in the
$(2n+4)$-polygon $P_n$ whose vertices are 
\[
 0,\;\alpha,\; \alpha\!+\!\beta_0,\;  \alpha\!+\!\beta_0\!+\!\beta_1,\; \dots\; ,
 \; \alpha\!+\!\beta_0\!+\!\dots\!+\!\beta_n, 
 \; \beta_0\!+\!\dots\!+\!\beta_n,\;  \beta_1\!+\!\dots\!+\!\beta_n,\;\dots\; ,\beta_n
\]
As in the previous example, the sets $\chi(\s_n,\Lambda \s_n)$ and $\chi(\s_n,\Lambda \f_n)$ coincide
(their multiplicities are related by a factor 4).
In the pictures below, the set of roots and the multiplicities of each weight in 
$\chi(\s_n,\Lambda \f_n)$ are described for the cases $n=4$ and $n=7$.

It follows from  \eqref{eqn:dixmier3} that
\[
H^p(\f_n)=(\Lambda^p V_n^*)^E \oplus (\Lambda^{p-1} V_n^*)_E\otimes E^*.
\]
This implies that the $\a$-isotopic component of $H^p(\f_n)$ of weight $\mu$, $H^p(\f_n)_{(\mu)}$, 
is generated by the set of $\gl(2)$-highest weight vectors of $\Lambda^p V_n^*$ of weight $\mu$
and by the set of $\gl(2)$-lowest weight vectors $\Lambda^{p-1} V_n^*$ of weight $\mu+(2,0)$.

The $\gl(2)$-module structure of $\Lambda^p V_n^*$ follows from 
the $\sl(2)$-character formula of $\Lambda^p V_n^*\simeq \Lambda^p V_n$ (see \cite[Ch 7]{S})
\[
 \text{ch}^p_n(q)=q^{-p(n+1-p)}\binom{[n+1]_{q^2}}{[p]_{q^2}},
\]
where $[a]_q=1+q+\dots+q^{a-1}$ is the $q$-analog of the number $a$.
The multiplicities of the highest weights are obtained by subtracting
two consecutive coefficients of $\text{ch}^p_n(q)$ but, as far as we know, 
there is no closed formula for these multiplicities.

In the pictures below, for $n=4$ and $n=7$, the negatives of the highest weights of the 
$\gl(2)$-module $\Lambda V_n^*$
are indicated with diamonds while the bullets indicate the negatives of the lowest weights
(translated by (2,0)).
The pictures actually show all the negatives of the $\a$-weights appearing in $H^*(\f_n)$, 
where diamonds correspond to 
(the negatives of) $\gl(2)$-highest weight vectors and bullets correspond to 
(the negatives of) $\gl(2)$-lowest weight vectors.
We point out that for $n=7$ there are few of these with multiplicity 2; 
this is also indicated in the picture.

The set $\{E,e_{n-2k}:\,k=0,\dots,n\}$ is a basis of $\{\f_n\}$.
Let $\{E^*,e_{n-2k}^*:\,k=0,\dots,n\}$ be its dual basis.
Among the highest and lowest weight vectors of $\Lambda^p V_n^*$ we have,
respectively, the following obvious ones:
\[
e_{2p-2-n}^*\wedge\dots\wedge e_{-n}^*\quad\text{ and }\quad e_{n}^*\wedge\dots\wedge e_{n+2-2p}^*,
\]
which yield the following two cohomology classes in $H^p(\f_n)$:
\[
e_{2p-2-n}^*\wedge\dots\wedge e_{-n}^*\quad\text{ and }\quad 
e_{n}^*\wedge\dots\wedge e_{n+4-2p}^*\wedge E^*,
\]
with respective weights $-\beta_{n-p+1}-\dots-\beta_{n}$ and $-\alpha-\beta_{0}-\dots-\beta_{p-2}$.

Finally, by Proposition \ref{prop:hoch-serre-abelian}, we have that 
\[
 H^p( \s_n,\lambda)=H^p( \f_n)_{(-\lambda)}\otimes\Lambda^0\a^*\;\oplus\;
 H^{p-1}(\f_n)_{(-\lambda)}\otimes\Lambda^1\a^*\;\oplus\;
 H^{p-2}( \f_n)_{(-\lambda)}\otimes\Lambda^2\a^*
\]
and thus $\lambda\in\Gamma_\circ^p$ if and only if $-\lambda$ is an 
$\mathfrak{a}$-weight of either $H^p(\f_n)$, $H^{p-1}(\f_n)$ or $H^{p-2}(\f_n)$. 
Therefore, $\Gamma_\circ(\s_n)$ is the set indicated with diamonds and bullets
in the pictures below for $n=4$ and $n=7$. 
In particular, all the $\s_n$-weights $\lambda$ that are vertices of the polygon 
$P_n=\chi(\s_n,\Lambda \s_n)$ are indeed in $\Gamma_\circ(\s_n)$.
\begin{center}

\tiny
 
\definecolor{zzttqq}{rgb}{0.6,0.2,0}
\definecolor{qqqqff}{rgb}{0,0,1}
\definecolor{cqcqcq}{rgb}{0.75,0.75,0.75}
\begin{tikzpicture}[line cap=round,line join=round,>=triangle 45,x=0.8cm,y=0.8cm]
\draw [color=cqcqcq,dash pattern=on 1pt off 1pt, xstep=0.8cm,ystep=0.8cm] (-1,-3) grid (8,4);
\clip(-1,-3) rectangle (8,4);
\fill[color=zzttqq,fill=zzttqq,fill opacity=0.1] (3,-2) -- (1,-1) -- (0,0) -- (0,1) -- (1,2) -- (3,3) -- (4,3) -- (6,2) -- (7,1) -- (7,0) -- (6,-1) -- (4,-2) -- cycle;
\draw [->] (3,-2) -- (4,-2);
\draw [->] (3,-2) -- (1,-1);
\draw [->] (3,-2) -- (2,-1);
\draw [->] (3,-2) -- (3,-1);
\draw [->] (3,-2) -- (4,-1);
\draw [->] (3,-2) -- (5,-1);
\draw [color=zzttqq] (3,-2)-- (1,-1);
\draw [color=zzttqq] (1,-1)-- (0,0);
\draw [color=zzttqq] (0,0)-- (0,1);
\draw [color=zzttqq] (0,1)-- (1,2);
\draw [color=zzttqq] (1,2)-- (3,3);
\draw [color=zzttqq] (3,3)-- (4,3);
\draw [color=zzttqq] (4,3)-- (6,2);
\draw [color=zzttqq] (6,2)-- (7,1);
\draw [color=zzttqq] (7,1)-- (7,0);
\draw [color=zzttqq] (7,0)-- (6,-1);
\draw [color=zzttqq] (6,-1)-- (4,-2);
\draw [color=zzttqq] (4,-2)-- (3,-2);
\draw (0.96,-0.58) node[anchor=north west] {1};
\draw (-0.06,1.42) node[anchor=north west] {1};
\draw (-0.04,0.38) node[anchor=north west] {1};
\draw (0.94,2.42) node[anchor=north west] {1};
\draw (2.96,3.39) node[anchor=north west] {1};
\draw (3.94,3.4) node[anchor=north west] {1};
\draw (5.93,2.38) node[anchor=north west] {1};
\draw (6.94,1.38) node[anchor=north west] {1};
\draw (6.94,0.38) node[anchor=north west] {1};
\draw (5.96,-0.58) node[anchor=north west] {1};
\draw (3.97,-1.63) node[anchor=north west] {1};
\draw (0.94,0.38) node[anchor=north west] {2};
\draw (0.94,1.4) node[anchor=north west] {2};
\draw (1.94,2.36) node[anchor=north west] {2};
\draw (2.94,2.36) node[anchor=north west] {2};
\draw (3.92,2.38) node[anchor=north west] {2};
\draw (4.93,2.33) node[anchor=north west] {2};
\draw (5.92,1.36) node[anchor=north west] {2};
\draw (5.93,0.34) node[anchor=north west] {2};
\draw (4.94,-0.66) node[anchor=north west] {2};
\draw (3.92,-0.64) node[anchor=north west] {2};
\draw (2.94,-0.63) node[anchor=north west] {2};
\draw (1.92,-0.63) node[anchor=north west] {2};
\draw (1.94,0.37) node[anchor=north west] {3};
\draw (1.92,1.4) node[anchor=north west] {3};
\draw (4.93,1.44) node[anchor=north west] {3};
\draw (4.93,0.38) node[anchor=north west] {3};
\draw (2.91,1.37) node[anchor=north west] {4};
\draw (3.92,1.37) node[anchor=north west] {4};
\draw (2.91,0.34) node[anchor=north west] {4};
\draw (3.92,0.34) node[anchor=north west] {4};
\draw (2.52,-2.13) node[anchor=north west] {(0,0)};
\draw (3.52,-2.13) node[anchor=north west] {(2,0)};
\draw (0.42,-1.14) node[anchor=north west] {(-4,1)};
\begin{scriptsize}
\fill [color=qqqqff] (3,-2) ++(-2.5pt,0 pt) -- ++(2.5pt,2.5pt)--++(2.5pt,-2.5pt)--++(-2.5pt,-2.5pt)--++(-2.5pt,2.5pt);
\fill [color=qqqqff] (4,-2) circle (2.5pt);
\fill [color=qqqqff] (1,-1) ++(-2.5pt,0 pt) -- ++(2.5pt,2.5pt)--++(2.5pt,-2.5pt)--++(-2.5pt,-2.5pt)--++(-2.5pt,2.5pt);
\fill [color=qqqqff] (0,0) ++(-2.5pt,0 pt) -- ++(2.5pt,2.5pt)--++(2.5pt,-2.5pt)--++(-2.5pt,-2.5pt)--++(-2.5pt,2.5pt);
\fill [color=qqqqff] (0,1) ++(-2.5pt,0 pt) -- ++(2.5pt,2.5pt)--++(2.5pt,-2.5pt)--++(-2.5pt,-2.5pt)--++(-2.5pt,2.5pt);
\fill [color=qqqqff] (1,2) ++(-2.5pt,0 pt) -- ++(2.5pt,2.5pt)--++(2.5pt,-2.5pt)--++(-2.5pt,-2.5pt)--++(-2.5pt,2.5pt);
\fill [color=qqqqff] (3,3) ++(-2.5pt,0 pt) -- ++(2.5pt,2.5pt)--++(2.5pt,-2.5pt)--++(-2.5pt,-2.5pt)--++(-2.5pt,2.5pt);
\fill [color=qqqqff] (4,3) circle (2.5pt);
\fill [color=qqqqff] (6,2) circle (2.5pt);
\fill [color=qqqqff] (7,1) circle (2.5pt);
\fill [color=qqqqff] (7,0) circle (2.5pt);
\fill [color=qqqqff] (6,-1) circle (2.5pt);
\fill [color=qqqqff] (2,0) ++(-2.5pt,0 pt) -- ++(2.5pt,2.5pt)--++(2.5pt,-2.5pt)--++(-2.5pt,-2.5pt)--++(-2.5pt,2.5pt);
\fill [color=qqqqff] (5,0) circle (2.5pt);
\fill [color=qqqqff] (5,1) circle (2.5pt);
\fill [color=qqqqff] (2,1) ++(-2.5pt,0 pt) -- ++(2.5pt,2.5pt)--++(2.5pt,-2.5pt)--++(-2.5pt,-2.5pt)--++(-2.5pt,2.5pt);
\end{scriptsize}
\end{tikzpicture}


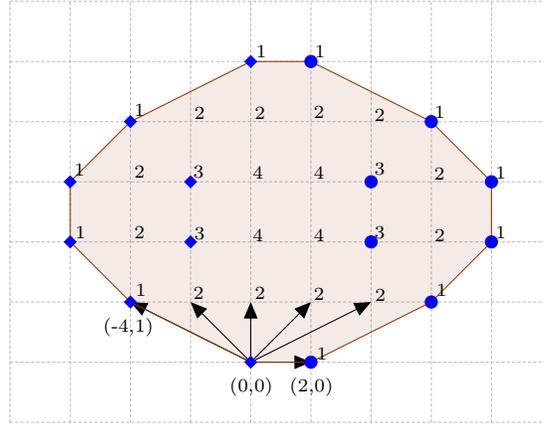
$\Gamma_\circ(\s_4)\subset P_4$ --- The negatives of the $\a$-weights in $H^*(\f_4)$ \\

\ \\


\tiny

\definecolor{qqqqff}{rgb}{0,0,1}
\definecolor{zzttqq}{rgb}{0.6,0.2,0}
\definecolor{cqcqcq}{rgb}{0.75,0.75,0.75}
\begin{tikzpicture}[line cap=round,line join=round,>=triangle 45,x=0.35cm,y=0.7cm]
\draw [color=cqcqcq,dash pattern=on 1pt off 1pt, xstep=0.35cm,ystep=0.7cm] (-16.5,-0.5) grid (18.5,8.5);
\clip(-16.5,-0.5) rectangle (18.5,8.5);
\fill[color=zzttqq,fill=zzttqq,fill opacity=0.1] (0,0) -- (2,0) -- (9,1) -- (14,2) -- (17,3) -- (18,4) -- (17,5) -- (14,6) -- (9,7) -- (2,8) -- (0,8) -- (-7,7) -- (-12,6) -- (-15,5) -- (-16,4) -- (-15,3) -- (-12,2) -- (-7,1) -- cycle;
\draw [color=zzttqq] (0,0)-- (2,0);
\draw [color=zzttqq] (2,0)-- (9,1);
\draw [color=zzttqq] (9,1)-- (14,2);
\draw [color=zzttqq] (14,2)-- (17,3);
\draw [color=zzttqq] (17,3)-- (18,4);
\draw [color=zzttqq] (18,4)-- (17,5);
\draw [color=zzttqq] (17,5)-- (14,6);
\draw [color=zzttqq] (14,6)-- (9,7);
\draw [color=zzttqq] (9,7)-- (2,8);
\draw [color=zzttqq] (2,8)-- (0,8);
\draw [color=zzttqq] (0,8)-- (-7,7);
\draw [color=zzttqq] (-7,7)-- (-12,6);
\draw [color=zzttqq] (-12,6)-- (-15,5);
\draw [color=zzttqq] (-15,5)-- (-16,4);
\draw [color=zzttqq] (-16,4)-- (-15,3);
\draw [color=zzttqq] (-15,3)-- (-12,2);
\draw [color=zzttqq] (-12,2)-- (-7,1);
\draw [color=zzttqq] (-7,1)-- (0,0);
\draw [->] (0,0) -- (-7,1);
\draw [->] (0,0) -- (-5,1);
\draw [->] (0,0) -- (-3,1);
\draw [->] (0,0) -- (-1,1);
\draw [->] (0,0) -- (1,1);
\draw [->] (0,0) -- (3,1);
\draw [->] (0,0) -- (5,1);
\draw [->] (0,0) -- (7,1);
\draw [->] (0,0) -- (2,0);
\draw (-7.25,1.45) node[anchor=north west] {1};
\draw (-5.25,1.45) node[anchor=north west] {2};
\draw (0.75,1.45) node[anchor=north west] {2};
\draw (2.75,1.45) node[anchor=north west] {2};
\draw (4.75,1.45) node[anchor=north west] {2};
\draw (6.75,1.45) node[anchor=north west] {2};
\draw (-1.25,1.45) node[anchor=north west] {2};
\draw (-3.25,1.45) node[anchor=north west] {2};
\draw (8.75,1.45) node[anchor=north west] {1};
\draw (11.75,2.45) node[anchor=north west] {2};
\draw (-12.25,2.45) node[anchor=north west] {1};
\draw (-10.25,2.45) node[anchor=north west] {2};
\draw (13.75,2.45) node[anchor=north west] {1};
\draw (-8.25,2.45) node[anchor=north west] {3};
\draw (-6.25,2.45) node[anchor=north west] {4};
\draw (-4.25,2.45) node[anchor=north west] {5};
\draw (-2.25,2.45) node[anchor=north west] {6};
\draw (-0.25,2.45) node[anchor=north west] {7};
\draw (1.75,2.45) node[anchor=north west] {7};
\draw (9.75,2.45) node[anchor=north west] {3};
\draw (7.75,2.45) node[anchor=north west] {4};
\draw (5.75,2.45) node[anchor=north west] {5};
\draw (3.75,2.45) node[anchor=north west] {6};
\draw (14.75,3.45) node[anchor=north west] {2};
\draw (12.75,3.45) node[anchor=north west] {3};
\draw (10.75,3.45) node[anchor=north west] {5};
\draw (8.75,3.45) node[anchor=north west] {7};
\draw (6.75,3.45) node[anchor=north west] {9};
\draw (4.75,3.45) node[anchor=north west] {11};
\draw (2.75,3.45) node[anchor=north west] {12};
\draw (0.75,3.45) node[anchor=north west] {12};
\draw (-1.25,3.45) node[anchor=north west] {12};
\draw (-3.25,3.45) node[anchor=north west] {11};
\draw (-5.25,3.45) node[anchor=north west] {9};
\draw (-7.25,3.45) node[anchor=north west] {7};
\draw (-9.25,3.45) node[anchor=north west] {5};
\draw (-11.25,3.45) node[anchor=north west] {3};
\draw (-15.25,3.45) node[anchor=north west] {1};
\draw (-13.25,3.45) node[anchor=north west] {2};
\draw (16.75,3.45) node[anchor=north west] {1};
\draw (16.75,3.45) node[anchor=north west] {1};
\draw (15.75,4.45) node[anchor=north west] {2};
\draw (13.75,4.45) node[anchor=north west] {3};
\draw (11.75,4.45) node[anchor=north west] {5};
\draw (9.75,4.45) node[anchor=north west] {8};
\draw (7.75,4.45) node[anchor=north west] {10};
\draw (5.75,4.45) node[anchor=north west] {12};
\draw (3.75,4.45) node[anchor=north west] {14};
\draw (1.75,4.45) node[anchor=north west] {15};
\draw (-14.25,4.45) node[anchor=north west] {2};
\draw (-12.25,4.45) node[anchor=north west] {3};
\draw (-10.25,4.45) node[anchor=north west] {5};
\draw (-8.25,4.45) node[anchor=north west] {8};
\draw (-6.25,4.45) node[anchor=north west] {10};
\draw (-4.25,4.45) node[anchor=north west] {12};
\draw (-2.25,4.45) node[anchor=north west] {14};
\draw (-0.25,4.45) node[anchor=north west] {15};
\draw (-16.25,4.45) node[anchor=north west] {1};
\draw (17.75,4.45) node[anchor=north west] {1};
\draw (14.75,5.45) node[anchor=north west] {2};
\draw (12.75,5.45) node[anchor=north west] {3};
\draw (10.75,5.45) node[anchor=north west] {5};
\draw (8.75,5.45) node[anchor=north west] {7};
\draw (6.75,5.45) node[anchor=north west] {9};
\draw (4.75,5.45) node[anchor=north west] {11};
\draw (2.75,5.45) node[anchor=north west] {12};
\draw (0.75,5.45) node[anchor=north west] {12};
\draw (-1.25,5.45) node[anchor=north west] {12};
\draw (-3.25,5.45) node[anchor=north west] {11};
\draw (-5.25,5.45) node[anchor=north west] {9};
\draw (-7.25,5.45) node[anchor=north west] {7};
\draw (-9.25,5.45) node[anchor=north west] {5};
\draw (-11.25,5.45) node[anchor=north west] {3};
\draw (16.75,5.45) node[anchor=north west] {1};
\draw (11.75,6.45) node[anchor=north west] {2};
\draw (-12.25,6.45) node[anchor=north west] {1};
\draw (-10.25,6.45) node[anchor=north west] {2};
\draw (13.75,6.45) node[anchor=north west] {1};
\draw (-8.25,6.45) node[anchor=north west] {3};
\draw (-6.25,6.45) node[anchor=north west] {4};
\draw (-4.25,6.45) node[anchor=north west] {5};
\draw (-2.25,6.45) node[anchor=north west] {6};
\draw (-0.25,6.45) node[anchor=north west] {7};
\draw (1.75,6.45) node[anchor=north west] {7};
\draw (9.75,6.45) node[anchor=north west] {3};
\draw (7.75,6.45) node[anchor=north west] {4};
\draw (5.75,6.45) node[anchor=north west] {5};
\draw (3.75,6.45) node[anchor=north west] {6};
\draw (-7.25,7.45) node[anchor=north west] {1};
\draw (-5.25,7.45) node[anchor=north west] {2};
\draw (0.75,7.45) node[anchor=north west] {2};
\draw (2.75,7.45) node[anchor=north west] {2};
\draw (4.75,7.45) node[anchor=north west] {2};
\draw (6.75,7.45) node[anchor=north west] {2};
\draw (-1.25,7.45) node[anchor=north west] {2};
\draw (-3.25,7.45) node[anchor=north west] {2};
\draw (8.75,7.45) node[anchor=north west] {1};
\draw (-15.25,5.45) node[anchor=north west] {1};
\draw (-13.25,5.45) node[anchor=north west] {2};
\draw (1.75,8.45) node[anchor=north west] {1};
\draw (-0.25,8.45) node[anchor=north west] {1};
\draw (-1.0,0.0) node[anchor=north west] {(0,0)};
\draw (1.0,-0.0) node[anchor=north west] {(2,0)};
\draw (-8.2,0.9) node[anchor=north west] {(-7,1)};
\begin{scriptsize}
\fill [color=qqqqff] (-7,1) ++(-2.0pt,0 pt) -- ++(2.0pt,2.0pt)--++(2.0pt,-2.0pt)--++(-2.0pt,-2.0pt)--++(-2.0pt,2.0pt);
\fill [color=qqqqff] (-12,2) ++(-2.0pt,0 pt) -- ++(2.0pt,2.0pt)--++(2.0pt,-2.0pt)--++(-2.0pt,-2.0pt)--++(-2.0pt,2.0pt);
\fill [color=qqqqff] (-15,3) ++(-2.0pt,0 pt) -- ++(2.0pt,2.0pt)--++(2.0pt,-2.0pt)--++(-2.0pt,-2.0pt)--++(-2.0pt,2.0pt);
\fill [color=qqqqff] (-16,4) ++(-2.0pt,0 pt) -- ++(2.0pt,2.0pt)--++(2.0pt,-2.0pt)--++(-2.0pt,-2.0pt)--++(-2.0pt,2.0pt);
\fill [color=qqqqff] (-15,5) ++(-2.0pt,0 pt) -- ++(2.0pt,2.0pt)--++(2.0pt,-2.0pt)--++(-2.0pt,-2.0pt)--++(-2.0pt,2.0pt);
\fill [color=qqqqff] (-12,6) ++(-2.0pt,0 pt) -- ++(2.0pt,2.0pt)--++(2.0pt,-2.0pt)--++(-2.0pt,-2.0pt)--++(-2.0pt,2.0pt);
\fill [color=qqqqff] (-7,7) ++(-2.0pt,0 pt) -- ++(2.0pt,2.0pt)--++(2.0pt,-2.0pt)--++(-2.0pt,-2.0pt)--++(-2.0pt,2.0pt);
\fill [color=qqqqff] (0,8) ++(-2.0pt,0 pt) -- ++(2.0pt,2.0pt)--++(2.0pt,-2.0pt)--++(-2.0pt,-2.0pt)--++(-2.0pt,2.0pt);
\fill [color=qqqqff] (-12,4) ++(-2.0pt,0 pt) -- ++(2.0pt,2.0pt)--++(2.0pt,-2.0pt)--++(-2.0pt,-2.0pt)--++(-2.0pt,2.0pt);
\fill [color=qqqqff] (-11,5) ++(-2.0pt,0 pt) -- ++(2.0pt,2.0pt)--++(2.0pt,-2.0pt)--++(-2.0pt,-2.0pt)--++(-2.0pt,2.0pt);
\fill [color=qqqqff] (-11,3) ++(-2.0pt,0 pt) -- ++(2.0pt,2.0pt)--++(2.0pt,-2.0pt)--++(-2.0pt,-2.0pt)--++(-2.0pt,2.0pt);
\fill [color=qqqqff] (-10,4) ++(-2.0pt,0 pt) -- ++(2.0pt,2.0pt)--++(2.0pt,-2.0pt)--++(-2.0pt,-2.0pt)--++(-2.0pt,2.0pt);
\fill [color=qqqqff] (-9,5) ++(-2.0pt,0 pt) -- ++(2.0pt,2.0pt)--++(2.0pt,-2.0pt)--++(-2.0pt,-2.0pt)--++(-2.0pt,2.0pt);
\fill [color=qqqqff] (-8,6) ++(-2.0pt,0 pt) -- ++(2.0pt,2.0pt)--++(2.0pt,-2.0pt)--++(-2.0pt,-2.0pt)--++(-2.0pt,2.0pt);
\fill [color=qqqqff] (-9,3) ++(-2.0pt,0 pt) -- ++(2.0pt,2.0pt)--++(2.0pt,-2.0pt)--++(-2.0pt,-2.0pt)--++(-2.0pt,2.0pt);
\fill [color=qqqqff] (-8,2) ++(-2.0pt,0 pt) -- ++(2.0pt,2.0pt)--++(2.0pt,-2.0pt)--++(-2.0pt,-2.0pt)--++(-2.0pt,2.0pt);
\fill [color=qqqqff] (-8,4) ++(-2.0pt,0 pt) -- ++(2.0pt,2.0pt)--++(2.0pt,-2.0pt)--++(-2.0pt,-2.0pt)--++(-2.0pt,2.0pt);
\fill [color=qqqqff] (-7,5) ++(-2.0pt,0 pt) -- ++(2.0pt,2.0pt)--++(2.0pt,-2.0pt)--++(-2.0pt,-2.0pt)--++(-2.0pt,2.0pt);
\fill [color=qqqqff] (-5,5) ++(-2.0pt,0 pt) -- ++(2.0pt,2.0pt)--++(2.0pt,-2.0pt)--++(-2.0pt,-2.0pt)--++(-2.0pt,2.0pt);
\fill [color=qqqqff] (-7,3) ++(-2.0pt,0 pt) -- ++(2.0pt,2.0pt)--++(2.0pt,-2.0pt)--++(-2.0pt,-2.0pt)--++(-2.0pt,2.0pt);
\fill [color=qqqqff] (-5,3) ++(-2.0pt,0 pt) -- ++(2.0pt,2.0pt)--++(2.0pt,-2.0pt)--++(-2.0pt,-2.0pt)--++(-2.0pt,2.0pt);
\fill [color=qqqqff] (-4,4) ++(-2.0pt,0 pt) -- ++(2.0pt,2.0pt)--++(2.0pt,-2.0pt)--++(-2.0pt,-2.0pt)--++(-2.0pt,2.0pt);
\fill [color=qqqqff] (-3,5) ++(-2.0pt,0 pt) -- ++(2.0pt,2.0pt)--++(2.0pt,-2.0pt)--++(-2.0pt,-2.0pt)--++(-2.0pt,2.0pt);
\fill [color=qqqqff] (-3,3) ++(-2.0pt,0 pt) -- ++(2.0pt,2.0pt)--++(2.0pt,-2.0pt)--++(-2.0pt,-2.0pt)--++(-2.0pt,2.0pt);
\fill [color=qqqqff] (-4,2) ++(-2.0pt,0 pt) -- ++(2.0pt,2.0pt)--++(2.0pt,-2.0pt)--++(-2.0pt,-2.0pt)--++(-2.0pt,2.0pt);
\fill [color=qqqqff] (0,2) ++(-2.0pt,0 pt) -- ++(2.0pt,2.0pt)--++(2.0pt,-2.0pt)--++(-2.0pt,-2.0pt)--++(-2.0pt,2.0pt);
\draw [color=qqqqff] (-8,4) ++(-4.0pt,0 pt) -- ++(4.0pt,4.0pt)--++(4.0pt,-4.0pt)--++(-4.0pt,-4.0pt)--++(-4.0pt,4.0pt);
\draw [color=qqqqff] (-4,4) ++(-4.0pt,0 pt) -- ++(4.0pt,4.0pt)--++(4.0pt,-4.0pt)--++(-4.0pt,-4.0pt)--++(-4.0pt,4.0pt);
\fill [color=qqqqff] (0,4) ++(-2.0pt,0 pt) -- ++(2.0pt,2.0pt)--++(2.0pt,-2.0pt)--++(-2.0pt,-2.0pt)--++(-2.0pt,2.0pt);
\fill [color=qqqqff] (0,6) ++(-2.0pt,0 pt) -- ++(2.0pt,2.0pt)--++(2.0pt,-2.0pt)--++(-2.0pt,-2.0pt)--++(-2.0pt,2.0pt);
\fill [color=qqqqff] (-4,6) ++(-2.0pt,0 pt) -- ++(2.0pt,2.0pt)--++(2.0pt,-2.0pt)--++(-2.0pt,-2.0pt)--++(-2.0pt,2.0pt);
\fill [color=qqqqff] (2,8) circle (2.0pt);
\fill [color=qqqqff] (9,7) circle (2.0pt);
\fill [color=qqqqff] (14,6) circle (2.0pt);
\fill [color=qqqqff] (10,6) circle (2.0pt);
\fill [color=qqqqff] (6,6) circle (2.0pt);
\fill [color=qqqqff] (2,6) circle (2.0pt);
\fill [color=qqqqff] (2,4) circle (2.0pt);
\fill [color=qqqqff] (5,5) circle (2.0pt);
\fill [color=qqqqff] (7,5) circle (2.0pt);
\fill [color=qqqqff] (9,5) circle (2.0pt);
\fill [color=qqqqff] (11,5) circle (2.0pt);
\fill [color=qqqqff] (13,5) circle (2.0pt);
\fill [color=qqqqff] (17,5) circle (2.0pt);
\fill [color=qqqqff] (18,4) circle (2.0pt);
\fill [color=qqqqff] (14,4) circle (2.0pt);
\fill [color=qqqqff] (12,4) circle (2.0pt);
\fill [color=qqqqff] (10,4) circle (2.0pt);
\fill [color=qqqqff] (6,4) circle (2.0pt);
\fill [color=qqqqff] (5,3) circle (2.0pt);
\fill [color=qqqqff] (7,3) circle (2.0pt);
\fill [color=qqqqff] (9,3) circle (2.0pt);
\fill [color=qqqqff] (11,3) circle (2.0pt);
\fill [color=qqqqff] (13,3) circle (2.0pt);
\fill [color=qqqqff] (14,2) circle (2.0pt);
\fill [color=qqqqff] (10,2) circle (2.0pt);
\fill [color=qqqqff] (6,2) circle (2.0pt);
\fill [color=qqqqff] (2,2) circle (2.0pt);
\fill [color=qqqqff] (9,1) circle (2.0pt);
\draw [color=qqqqff] (10,4) circle (3.5pt);
\draw [color=qqqqff] (6,4) circle (3.5pt);
\fill [color=qqqqff] (17,3) circle (2.0pt);
\fill [color=qqqqff] (2,0) circle (2.0pt);
\fill [color=qqqqff] (0,0) ++(-2.0pt,0 pt) -- ++(2.0pt,2.0pt)--++(2.0pt,-2.0pt)--++(-2.0pt,-2.0pt)--++(-2.0pt,2.0pt);
%
\end{scriptsize}
\end{tikzpicture}


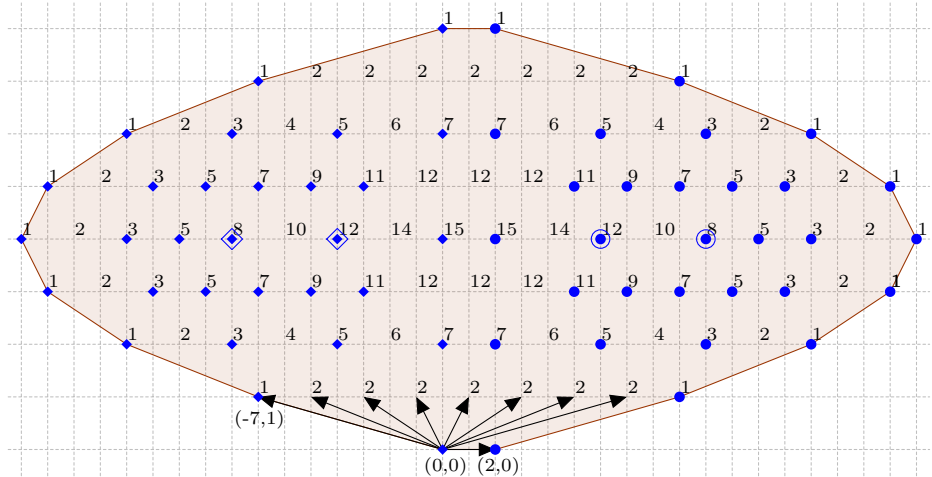
$\Gamma_\circ(\s_7)\subset P_7$ --- The negatives of the  $\a$-weights in $H^*(\f_7)$

\end{center}

\subsection{An approximation to $\Gamma_\circ^1(\s)$}
We now prove the following result.

\begin{theorem}
If $\s$ is any solvable Lie algebra and $\s''$ is 
the second term of the derived series of $\s$, then
\[ 
\Gamma_\circ^1(\s)\subset\chi(\s, \s/\s''). 
\]  
\end{theorem}

\begin{proof}
Let 
\[
 \s=\mathfrak{c}\oplus\sum_{\alpha\in\Phi(\s,\mathfrak{c})}\s_{\alpha}
\]
be the (generalized) root decomposition of $\s$ with respect to a 
Cartan subalgebra $\mathfrak{c}$ of $\s$, being $\Phi(\s,\mathfrak{c})$ 
the corresponding set of roots. 
It is clear that 
\begin{equation}\label{eq.s_s''1}
\s'=\big(\mathfrak{c}'+\mathfrak{c}_1\big)
\oplus\sum_{\alpha\in\Phi(\s,\mathfrak{c})}\s_{\alpha},
\end{equation}
where $\mathfrak{c}_1=  \sum_{\beta\in\Phi(\s,\mathfrak{c})}[\s_{\beta},\s_{-\beta}]$ 
is an ideal of $\mathfrak{c}$, and 
\begin{equation}\label{eq.s_s''2}
\s'' =\big(\mathfrak{c}''+\mathfrak{c}_1\big)\oplus
\Big(\sum_{\alpha\in\Phi(\s,\mathfrak{c})}[\mathfrak{c}'+\mathfrak{c}_1,\s_{\alpha}]
+\sum_{{\beta_1,\beta_2\in\Phi(\s,\mathfrak{c})} \atop {\beta_1+\beta_2\ne0}}
[\s_{\beta_1},\s_{\beta_2}]\Big).
\end{equation}
Let $\lambda\in\Gamma_\circ^1(\s)$.
Since the weight zero is in $\chi(\s, \s/\s'')$ we may assume that $\lambda\ne0$.
We know, from Theorem \ref{thm:TH}, that 
there exists $\alpha_0\in\Phi(\s,\mathfrak{c})$ 
such that $\lambda|_{\mathfrak{c}}=\alpha_0$.
Let $\h_0=\ker \alpha_0\subsetneq\mathfrak{c}$ and $D\in \mathfrak{c}$ so that $\alpha_0(D)\ne0$.
Therefore $\s=D\ltimes \h$, with  $\h=\ker\lambda$, and clearly 
\begin{align}\label{eq.h_h'1}
  \h &=\mathfrak{h}_0\oplus\sum_{\alpha\in\Phi(\s,\mathfrak{c})}\s_{\alpha}, \\\label{eq.h_h'2}
  \h' &=\mathfrak{h}_0'\oplus\sum_{\alpha\in\Phi(\s,\mathfrak{c})}[\h_0,\s_{\alpha}]
+\sum_{\beta_1,\beta_2\in\Phi(\s,\mathfrak{c})}[\s_{\beta_1},\s_{\beta_2}].
\end{align}
We know from \eqref{eqn:dixmier2} that 
$H^1(\s,\lambda)\simeq H^1(\h)^{\theta_0(D)+\alpha_0(D)} + H^0(\h)_{\theta_0(D)+\alpha_0(D)}$.
Since $\lambda\ne 0$ ($\alpha_0(D)\ne 0$) it follows that $H^0(\h)_{\theta_0(D)+\alpha_0(D)}=0$ and hence
\[ 
H^1(\s,\lambda)=H^1(\h)^{\theta_0(D)+\alpha_0(D)}. 
\]
Since $H^1(\h)=\big(\h/\h'\big)^*$ and $\lambda\in\Gamma_\circ^1(\s)$, we obtain that
$\alpha_0(D)$ is an eigenvalue of $D$ on $\h/\h'$.
It follows from \eqref{eq.h_h'1} and \eqref{eq.h_h'2} 
that $\alpha_0(D)=\gamma(D)$ for some 
$\gamma\in\Phi(\s,\mathfrak{c})$ such that 
\begin{equation}\label{eq.s_gamma}
  \s_{\gamma}\supsetneq [\h_0,\s_{\gamma}] + 
\sum_{{\beta_1,\beta_2\in\Phi(\s,\mathfrak{c})} \atop {\beta_1+\beta_2=\gamma}}
[\s_{\beta_1},\s_{\beta_2}].
\end{equation}
This fact implies that $\gamma|_{\h_0}=0$.
Since $\h_0=\ker\alpha_0$ and $\alpha_0(D)=\gamma(D)$,
it follows that $\gamma=\alpha_0$.
Finally, since $\mathfrak{c}'+\mathfrak{c}_1\subset\mathfrak{s}'\cap \mathfrak{c}\subset\h_0$, and hence
\[
[\mathfrak{c}'+\mathfrak{c}_1,\s_{\alpha_0}] \subset [\h_0,\s_{\alpha_0}],
\]
it follows from \eqref{eq.s_gamma} that 
\[
 \s_{\alpha_0}\supsetneq [\mathfrak{c}'+\mathfrak{c}_1,\s_{\alpha_0}] + 
\sum_{{\beta_1,\beta_2\in\Phi(\s,\mathfrak{c})} \atop {\beta_1+\beta_2=\alpha_0}}
[\s_{\beta_1},\s_{\beta_2}].
\]
This fact, together with \eqref{eq.s_s''2}, implies that 
$\alpha_0$ is a $\mathfrak{c}$-weight of $\s/s''$ and, in other words,
$\lambda\in\chi(\s, \s/\s'')$. 
\end{proof}

\section{Total cohomology}\label{sec:total}

Theorem \ref{thm:TH} allows us to define the \emph{Total cohomology} of $\g$ as
\[ 
TH^p(\g)=\bigoplus_{(\pi,V)\in\Gamma_\circ^p(\g)}H^p(\g,V).
\]
We know, from the results of Whitehead and Dixmier, that 
$TH^p(\g)=H^p(\g)$ if $\g$ is semisimple or nilpotent.

Let $\g=\g_0\ltimes \s$ be a Levi decomposition of a Lie algebra $\g$.
We will show now that it follows from Theorem \ref{thm:TH} that
\[
\text{Der}(\g)_0=\{D\in\text{Der}(\g):D(\g_0)\subseteq\g_0\}
\]
acts on $H^p(\g,V)$ for any $(\pi,V)\in\text{Irrep}(\g)$
and in particular the semisimple factor of a Levi decomposition of $\text{Der}(\g)$ does. 

In fact, if $(\pi,V)\not\in \Gamma_\circ(\g)$ then $H^*(\g,V)=0$ and thus
we may assume that $(\pi,V)\in\Gamma_\circ^p(\g)$ for some $p$.
Given $D\in\text{Der}(\g)_0$, let $X_D\in \g_0$ such that 
$D|_{\g_0}=\text{ad}_{\g_0}(X_D)$. 
It is clear that that $D\mapsto X_D$ is a well defined Lie algebra homomorphism.

Let $\text{Der}(\g)_0$ act on $\Lambda\g^*\otimes V$ via $\psi$ as follows:
if $D\in\text{Der}(\g)_0$, $v\in V$ and $f\in\Lambda\g^*$, then 
\begin{equation}\label{eq.def_theta_0}
 \psi(D)(f\otimes v) = - \overline{D}^{t}(f)\otimes v+f\otimes \pi(X_D)(v),
\end{equation}
where $\overline{D}$ is the canonical extension of $D$ to the exterior algebra of $\g$.
Notice that if $\g=\s$ is solvable, then $\psi=\theta_0$ (see \S\ref{sec:preliminaries}).

We need to show that $\psi(D)$ commutes with the coboundary operator on $\Lambda\g^*\otimes V$.
Since $D$ is a derivation, this is equivalent to showing that 
\begin{equation}\label{eq.commute}
 [\psi(D),\pi(X)]=\pi(D(X))
 \end{equation}
on $\Lambda^0\g\otimes V$, for all $D\in\text{Der}(\g)_0$ and $X\in\g$. 
If $X=X_{\g_0}+X_{\s}$ with $X_{\g_0}\in\g_0$ and $X_{\s}\in\s$ then, on $\Lambda^0\g\otimes V$,
$\psi(D)=\pi(X_D)$ and 
\begin{align*}
D(X)
& =\text{ad}_{\g_0}(X_D)(X_{\g_0})+D(X_{\s}) \\
& =[X_D,X-X_{\s}]+D(X_{\s}) \\
& =[X_D,X]+\tilde D(X_{\s})
\end{align*}
where $\tilde D=D-\text{ad}_{\g}(X_D)$.
Therefore, \eqref{eq.commute} is equivalent to 
\[
 [\pi(X_D),\pi(X)]=\pi([X_D,X]+\tilde D(X_{\s})).
\]
Since $\pi$ is a Lie algebra homomorphism, this is equivalent to showing
that $\pi(\tilde D(X_{\s}))=0$. 
Since $\pi\in\Gamma_\circ^p(\g)\subset\chi(\g,\Lambda\s)$,
this follows from the following lemma by using a Jordan-H\"older series of $\pi|_\s$.
\begin{lemma}\label{lemma.D_trivial}
If $\s$ is a solvable Lie algebra and $\sigma\in\chi(\s,\Lambda\s)$, 
then $\sigma\circ D=0$ for all $D\in\text{Der}(\s)$.
\end{lemma}
\begin{proof}
We know that $D(\s)$ is contained in the nilradical of 
$\s$ for  all $D\in\text{Der}(\s)$
(see \cite[Corollary 2, Ch II.7]{J}) and thus $\text{ad}(D(X))$ is nilpotent.
Therefore $\sigma(D(X))=0$ for all $\sigma\in\chi(\s,\Lambda\s)$.
\end{proof}

We have then proved the following theorem.
\begin{theorem}
Let $\g=\g_0\ltimes\s$ be the Levi decomposition of a 
Lie algebra $\g$ and let 
$\text{Der}(\g)_0=\{D\in\text{Der}(\g):D(\g_0)\subseteq\g_0\}$.
If $(\pi,V)\in\Gamma_\circ^p(\g)$, then 
the action $\psi$ of $\text{Der}(\g)_0$ on  $\Lambda\g^*\otimes V$ 
defined in $\eqref{eq.def_theta_0}$ induces an action of 
$\text{Der}(\g)_0$ on $H^p(\g,V)$ for all $p$.
In particular  $\psi$ defines an action of 
$\text{Der}(\g)_0$ on  $TH^p(\g)$ for all $p$. 
\end{theorem}

\section{Total cohomology of solvable Lie algebras}

Let $\s$ be a solvable Lie algebra.
The set of characters of $\s$, $\HomL(\s,\mathbb K)$, 
can be identified with the dual space $(\s/\s')^*$
and it follows from Theorem \ref{thm:TH} that 
\[
 \Gamma_\circ^p(\s)\subseteq\chi(\s,\Lambda^p\s)
=\Big\{\textstyle\sum_{i=1}^p\lambda_i:
\lambda_i\in\chi(\s,\s)\Big\}\subseteq\HomL(\s,\mathbb K)
\]
for all $p=0,\dots,\dim\s$ (cf.\ \cite[\S 3]{M}). 

Let $\h$ be an ideal of $\s$ of codimension 1 and let $D\not\in\h$ so that $\s=D\ltimes\h$.
For any character $\lambda$ of $\s$, $\lambda\in \HomL(\h,\K)$, let $\theta_0(D)$ and $\theta_\lambda(D)$
be defined as in \S \ref{sec:preliminaries} and  
consider, for each $p=0,\dots,\dim\h$, the subset $R^p\subseteq \K\oplus\HomL(\h,\K)$ defined by
\[
R^p=
\{(z,\lambda): \lambda\in\Gamma_\circ^{p}(\h)
\text{ and $z$ is an eigenvalue of $-\theta_0(D)$ in $H^{p}(\h,\lambda)$}\}.
\]

\begin{lemma}
Let $\s$ be a solvable Lie algebra, $\h$ an ideal of codimension 1 and let $D\not\in\h$.
Then, under the map 
\begin{equation}\label{def.identification}
\begin{split} 
\HomL(\s,\K)&\rightarrow \K\oplus\HomL(\h,\K) \\
\lambda\;\;&\mapsto\,\;\big(\lambda(X), \lambda|_\h\big)
\end{split}
\end{equation}
the sets $\Gamma_\circ^p(\s)$ and $R^{p}\cup R^{p-1}$
are in bijective correspondence.
\end{lemma}

\begin{proof}
By \eqref{eqn:dixmier2} we have that
\[
H^p(\s,\lambda) \simeq 
H^{p}(\h,\lambda)^{\theta_0(D)+\lambda(D)} \oplus
H^{p-1}(\h,\lambda)_{\theta_0(D)+\lambda(D)}.
\]
From this it follows that if $\lambda\in \Gamma_\circ^p(\s)$, then 
$\big(\lambda(D), \lambda|_\h\big)$ belongs to $R^{p}$ or $R^{p-1}$.

Conversely, given $(z,\mu)\in R^p$ 
let $\lambda$ be define by $\lambda(D)=z$ and $\lambda|_\h=\mu$.
It follows that $\lambda\in \Gamma_\circ^p$ (and also $\lambda\in \Gamma_\circ^{p+1}$ if $p\le\dim \h$),
provided that $\lambda\in\HomL(\s,\K)$.
Indeed, $\s'=D(\h)+\h'$ and since $\mu|_{\h'}=0$ and $\mu\circ D=0$ (see Lemma \ref{lemma.D_trivial}),
we obtain that $\mu|_{\s'}=0$ as needed.
\end{proof}

We prove now one of the main results of this paper
stating that eliminating (modifying) the semisimple part of
$\ad D$ in a solvable Lie algebra $D\ltimes\h$, does not affect the total cohomology.
This phenomenon will later on the paper be understood in the framework of linear deformations
of Lie algebras.

Recall that the semisimple and nilpotent parts of any derivation of an arbitrary Lie algebra
are also derivations (see for instance \cite[Lemma 4.2B]{H}).

\begin{theorem}\label{thm:induction}
Let $\s$ be a solvable Lie algebra and let $\h$ be an ideal of codimension 1. 
Given $D\not\in\h$, let $D=D_s+D_n$ be the Jordan decomposition of $\ad D$ in $\Der(\h)$
and let $\s_n=D_n\ltimes\h$. 
Then
\[
\dim TH^p(\s)=\dim TH^p(\s_n)
\]
for all $p=0,\dots,\dim\s$.
\end{theorem}

\begin{proof}
From the previous lemma (and its proof) it follows that 
\begin{align*}
TH^p(\s)
&=\bigoplus_{\lambda\in\Gamma_\circ^p(\s)}H^p(\s,\lambda) \\
&\simeq
\bigoplus_{(z,\mu)\in R^{p}}H^{p}(\h,\mu)^{\theta_0(D)+z}
\bigoplus_{(z,\mu)\in R^{p-1}}H^{p-1}(\h,\mu)_{\theta_0(D)+z}.
\end{align*} 

Since $\theta_0(D)=\theta_0(D_s)+\theta_0(D_n)$ is the Jordan decomposition of $\theta_0(D)$ on
$H^*(\h,\mu)$ it follows that 
\begin{align*}
\bigoplus_{(z,\mu)\in R^{p}}H^{p}(\h,\mu)^{\theta_0(D)+z}
& = \bigoplus_{\mu\in\Gamma_\circ^{p}(\h)}H^{p}(\h,\mu)^{\theta_0(D_n)} \\ 
\bigoplus_{(z,\mu)\in R^{p-1}}H^{p-1}(\h,\mu)_{\theta_0(D)+z} 
& = \bigoplus_{\mu\in\Gamma_\circ^{p-1}(\h)}H^{p-1}(\h,\mu)_{\theta_0(D_n)} .
\end{align*} 
Therefore 
\begin{align*}
TH^p(\s)
& \simeq\bigoplus_{\mu\in\Gamma_\circ^{p}(\h)}H^{p}(\h,\mu)^{\theta_0(D_n)}
\bigoplus_{\mu\in\Gamma_\circ^{p-1}(\h)}H^{p-1}(\h,\mu)_{\theta_0(D_n)} \\
&\simeq\bigoplus_{\lambda\in\Gamma^p(\s_n)}H^p(\s_n,\lambda) \\
&=TH^p(\s_n)
\end{align*}
and the proof is complete. 
\end{proof}

\begin{corollary}\label{coro.Total_Ssimple}
Let $\h$ be a solvable Lie algebra and $D\in\Der(\h)$.
Then, for any semisimple derivation 
$S\in\text{Cent}_{\Der(\h)}(D)$, 
\begin{equation*}
 \dim TH^p((S+D)\ltimes \h)=\dim TH^p(D\ltimes \h),
\end{equation*}
for all $p=0,\dots,\dim(\h)+1$.
\end{corollary}
\begin{proof}
Since $(S+D)_n=D_n$, both sides are equal to $\dim TH^p(D_n\ltimes \h)$.
\end{proof}

If $\s$ is a solvable Lie algebra which is not nilpotent, then $\s$ is 
the semidirect product $D\ltimes\h$ where $\h$ is an ideal of $\s$
(of codimension 1) containing the nilradical of $\s$ and $D$ is a non-nilpotent derivation of $\h$.
Now, Theorem \ref{thm:induction} implies that $\dim TH^p(\s)=\dim TH^p(D_n\ltimes\h)$ for all $p$, 
where $D_n$ is the nilpotent part of $D$.
Since the nilradical of $D_n\ltimes\h$ is larger than the nilradical of $\s$,
an inductive argument yields a nilpotent Lie algebra $\n$ (of the same dimension of $\s$)
such that
\[ \dim TH^p(\s) = \dim H^p(\n), \]
for all $p=0,\dots,\dim\s$.
In \S\ref{sec:nilshadow} we will prove that this nilpotent Lie algebra $\n$ is the nilshadow of $\s$.

\section{Linear deformations}\label{sec:deformations}

In this section we study some aspects of the deformation theory of Lie algebras
and thus we will think of a Lie algebra structure on a given 
vector space $V$ over $\K$ (non-necessarily algebraically closed)
as a bilinear skew-symmetric map on $V$
\[ \mu: \Lambda^2 V \longrightarrow V, \]
satisfying the Jacobi identity
\[  \circlearrowleft\mu(\mu(x,y),z)=0, \]
for all $x,y,z\in V$ (the symbol $\circlearrowleft$ means cyclic sum).
Notice that $\mu$ is a Lie algebra if and only if $t\mu$ is a Lie algebra for all $t\in\K$.

Given a  Lie algebra $\mu$ on $V$, 
a family of Lie algebra structures on $V$
\[ \mu_t=\mu+t\sigma,\quad t\in\K,\]
where $\sigma:\Lambda^2 V \longrightarrow V$,
will be called a \emph{linear deformation} of $\mu$.
We shall also say that the members of the family are linear deformations of $\mu$.

Since $\mu$ satifies the Jacobi identity, it follows that the condition
\[ \circlearrowleft\mu_t(\mu_t(x,y),z)=0,\quad \text{for all } t\in\K,\]
is equivalent to both
\[ \circlearrowleft\sigma(\sigma(x,y),z)=0, \]
and
\begin{equation}\label{eqn:2-cocycle}
 \circlearrowleft\mu(\sigma(x,y),z)+\circlearrowleft\sigma(\mu(x,y),z))=0.
\end{equation}
That is, $\sigma$ is a Lie algebra structure which is a 2-cocycle of $\mu$.
In this case we shall say that $\sigma$ is an \emph{infinitesimal linear deformation} of $\mu$.
Notice that condition (\ref{eqn:2-cocycle}) is symmetric in $\mu$ and $\sigma$ and 
hence $\mu$ is a 2-cocycle of $\sigma$ if and only if $\sigma$ is a 2-cocycle of $\mu$
and we have the following proposition.

\begin{proposition}
 Given two Lie algebra structures on $V$, $\mu$ and $\sigma$, the following are equivalent.
\begin{enumerate}[(1)]
 \item $\sigma$ is a 2-cocycle of $\mu$.
 \item $\mu$ is a 2-cocycle of $\sigma$.
 \item $\mu_t=\mu+t\sigma$ is a Lie algebra structure for all $t\in\K$.
 \item $\sigma_t=\sigma+t\mu$ is a Lie algebra structure for all $t\in\K$.
 \item $\mu_1=\mu+\sigma$ is a Lie algebra structure.
 \item $\sigma_1=\sigma+\mu$ is a Lie algebra structure.
 \item $t_1\mu+t_2\sigma$ is a Lie algebra structure for all $t_1,t_2\in\K$.
\end{enumerate}
\end{proposition}

\subsection{Lie algebra extensions}\label{subsec:extensions}

Given two Lie algebras $\a$ and $\h$, the extensions of $\a$ by $\h$
\[
 0\;{\longrightarrow}\;\h\;\overset{\iota}{\longrightarrow}\;
\g\;\overset{\pi}{\longrightarrow}\;\a\;\longrightarrow\; 0
\]
are in correspondence with pairs of linear maps 
\[ \alpha:\a \to \Der(\h) \quad\text{and}\quad  \rho :\Lambda^2\a \to \h \]
satisfying, for all $X,Y,Z\in\a$,
\begin{eqnarray} 
[\alpha(X),\alpha(Y)]_{\End{\h}} &=& \alpha([X,Y]_{\a})+\ad_{\h}(\rho(X\wedge Y)),\label{eqn:extensions1}\\ 
\circlearrowleft\alpha(X)\big(\rho(Y\wedge Z)\big) &=& \circlearrowleft\rho([X,Y]_{\a}\wedge Z).\label{eqn:extensions2}
\end{eqnarray}
Given such an extension and a (linear) section $s:\a\to\g$, the linear maps
\begin{align*}
&\alpha:\a \to \Der(\h), & &\alpha(X)=\ad(s(X))|_{\h} \\
&\rho :\Lambda^2\a \to \h, & &\rho(X\wedge Y)=[s(X),s(Y)]-s([X,Y])
\end{align*}
satisfy (\ref{eqn:extensions1}) and (\ref{eqn:extensions2}).
Conversely, given maps $\alpha$ and $\rho$ satisfying 
(\ref{eqn:extensions1}) and (\ref{eqn:extensions2}), they define on the vector space $\g=\a\oplus\h$
a Lie algebra structure, that we will denote by 
\[
 \g_{\alpha,\rho}=\a\oplus_{\alpha,\rho}\h,
\]
defined by: 
\begin{enumerate}[(1)]
 \item $[X,Y]_{\g}=[X,Y]_{\a} + \rho(X\wedge Y)$, for $X,Y\in\a$,
 \item $[X,H]_{\g}=\alpha(X)(H)$, for $X\in\a$, $H\in\h$,
 \item $[H,K]_{\g}=[H,K]_{\h}$, for $H,K\in\h$.
\end{enumerate}
This new Lie bracket on $\a\oplus\h$ coincides on the ideal $\h$ with the original one, and in turn it induces
on the quotient $\a$ the original one also.
 
We point out that in the particular case when $\rho=0$, it turns out that
$\alpha$ is a Lie algebra homomorphism defining an action of $\a$ on $\h$,
and $\g_{\alpha,\rho}=\a\ltimes_\alpha\h$
is the semidirect product corresponding to this action.

\subsection{Linear deformations of extensions and their cohomology}\label{subsec:linear-def-ext}

Let $\a$ and $\h$ be two given Lie algebras.
Let $\g_{\alpha_1,\rho_1}$ and $\g_{\alpha_2,\rho_2}$ be two extensions of $\a$ by $\h$ 
and $[\cdot,\cdot]_1$ and $[\cdot,\cdot]_2$ be the corresponding Lie brackets.

It follows from \S\ref{subsec:extensions} that 
$[\cdot,\cdot]_2-[\cdot,\cdot]_1$ is a Lie bracket on $\a\oplus\h$ 
if and only if the linear maps
\begin{align*}
 \beta&=\alpha_2-\alpha_1:\a \to \Der(\h) \\
 \tau&=\rho_2-\rho_1 :\Lambda^2\a \to \h
\end{align*}
define an extension $\a_{\text{ab}}\oplus_{\beta,\tau}\h_{\text{ab}}$ of $\a_{\text{ab}}$ by $\h_{\text{ab}}$,
the vector spaces $\a$ and $\h$ with abelian Lie algebra structures.
This is, according to (\ref{eqn:extensions1}) and (\ref{eqn:extensions2}), equivalent to 
\begin{equation}\label{eq.dataAb}
\begin{split} 
&[\beta(X),\beta(Y)]_{\End{\h}}=0,\\
&\circlearrowleft\beta(X)\big(\tau(Y\wedge Z)\big)=0,
\end{split}
\end{equation}
for all $X,Y,Z\in\a$.
Translating condition (\ref{eqn:extensions1}) for $\alpha_2=\alpha_1+\beta$ and 
condition (\ref{eqn:extensions2}) for $\rho_2=\rho_1+\tau$ in terms of $\alpha_1$, $\rho_1$, $\beta$ and 
$\tau$ we obtain the following proposition
that states conditions on $\beta$ and $\tau$ such that 
the extension $\a_{\text{ab}}\oplus_{\beta,\tau}\h_{\text{ab}}$
is an infinitesimal linear deformation of $\a\oplus_{\alpha,\sigma}\h$.

\begin{proposition}Given two Lie algebras $\a$ by $\h$, 
let $\a\oplus_{\alpha,\sigma}\h$ be an extension of $\a$ by $\h$ and let 
\[ \beta:\a \to \Der(\h) \quad\text{and}\quad \tau :\Lambda^2\a \to \h \]
be given linear maps.
Then  $\alpha+t\beta$ and $\rho+t\tau$ define an extension
of $\a$ by $\h$ for all $t\in\K$ if and only if the following conditions are satisfied for all $X,Y,Z\in\a$:
\begin{enumerate}[(1)]
\item $[\beta(X),\beta(Y)]_{\End{\h}}=0$, 
\item $\circlearrowleft\beta(X)\big(\tau(Y\wedge Z)\big)=0$,
\item $[\alpha(X),\beta(Y)]_{\End{\h}}+[\beta(X),\alpha(Y)]_{\End{\h}}=
\beta([X,Y]_{\a})+\ad_{\h}(\tau(X\wedge Y))$,
\item $\circlearrowleft\beta(X)\big(\rho(Y\wedge Z)\big)\;+
\circlearrowleft\alpha(X)\big(\tau(Y\wedge Z)\big)
=\circlearrowleft\tau([X,Y]_{\a}\wedge Z)$.
\end{enumerate}
\end{proposition}

\begin{remark} The proof is straightforward, notice that 
conditions (1) and (2) say that  $\beta$ and $\tau$ define an extension of $\a_{\text{ab}}$ by $\h_{\text{ab}}$ 
 and conditions (3) and (4) assure that $\alpha+\beta$ and $\rho+\tau$ define an extension of $\a$ by $\h$.
\end{remark}

We will be particularly interested in the case when $\a=\a_{\text{ab}}$ is abelian and $\tau=0$.
In other words, we want to look at infinitesimal linear deformations of the form 
$\a\ltimes_\beta\h_{\text{ab}}$ of 
an extension of an abelian Lie algebra $\a$ by $\h$.
The above proposition says that $\beta:\a\to\Der(\h)$ must satisfy the following 
properties for $X,Y,Z\in\a$:
\begin{enumerate}[(1')]
\item $[\beta(X),\beta(Y)]_{\End{\h}}=0$, i.e., $\beta$ is a Lie algebra homomorphism and 
thus $\a$ acts on $\h$ by derivations.
\item (empty),
\item $[\alpha(X),\beta(Y)]_{\End{\h}}+[\beta(X),\alpha(Y)]_{\End{\h}}=0$,
\item ${\circlearrowleft}\beta(X)\big(\rho(Y\wedge Z)\big)=0$.
\end{enumerate}

\begin{remark}
 We point out that in case $\dim\a=1$, all conditions (1)-(4) are trivially fulfilled 
since $\a$ is abelian and $\tau=0$. 
These `codimension 1' linear deformations are the main subject studied in \cite{GO}.
We can think of deformations satisfying conditions (1')-(4') as generalizations
of those just mentioned.
\end{remark}

We consider now an even more particular type of linear deformations
of a given  extension.
\begin{definition}
Let $\h$ be a Lie algebra, $\a$ an abelian Lie algebra and 
$\g_{\alpha,\rho}=\a\oplus_{\alpha,\rho}\h$ an extension of $\a$ by $\h$ as in \S\ref{subsec:extensions}.
We say that the infinitesimal linear deformation $\a\oplus_{\beta,0}\h_{\text{ab}}$ 
of $\g_{\alpha,\rho}$ is \emph{elementary} if, for all $X,Y,Z\in\a$, it holds:
\begin{enumerate}[(3'')]
\item[(3'')]\label{(3'')} $[\alpha(X),\beta(Y)]_{\End{\h}}=0$,
\item[(4'')]\label{(4'')} $\beta(X)\big(\rho(Y\wedge Z)\big)=0$.
\end{enumerate}
Notice that necessarily, $\beta:\a \to \Der(\h)$
is a Lie algebra homomorphism.
We will say that $\g_{\alpha+\beta,\rho}$ is an \emph{elementary linear deformation} of $\g_{\alpha,\rho}$.
\end{definition}
\begin{proposition}\label{prop.Elementary}
Let $\h$ be a Lie algebra, $\a$ an abelian Lie algebra, 
$\g_{\alpha,\rho}$ an extension of $\a$ by $\h$ and $\a\oplus_{\beta,0}\h_{\text{ab}}$ 
an elementary infinitesimal linear deformation of 
$\g_{\alpha,\rho}$.
Let $\a=\a_1\oplus\a_2$ be a vector space decomposition of $\a$ and let
\[
 \beta_i(X)=
\begin{cases}
 \beta(X), & \text{if }X\in\a_i; \\
 0,        & \text{if }X\in\a_j,\; j\ne i.
\end{cases}
\]
(Clearly $\beta=\beta_1+\beta_2$.) 
Then $\a\oplus_{\beta_1,0}\h_{\text{ab}}$  is an 
elementary infinitesimal linear deformation of $\g_{\alpha,\rho}$ and   
$\a\oplus_{\beta_2,0}\h_{\text{ab}}$ is an 
elementary infinitesimal linear deformation of $\g_{\alpha+\beta_1,\rho}$.
In particular, $\g_{t\alpha+t_1\beta_1+t_2\beta_2,t\rho}$ is a Lie algebra for all
$t,t_1,t_2\in\K$.
\end{proposition}

\begin{proof} It is clear that both $\beta_1$ and $\beta_2$ are Lie algebra homomorphisms.
Since the following identities hold 
\begin{itemize}
 \item[] $[\alpha(X),\beta_1(Y)]_{\End{\h}}=0$,
 \item[] $[\alpha(X)+\beta_1(X),\beta_2(Y)]_{\End{\h}}=0$,
 \item[] $\beta_i(X)\big(\rho(Y\wedge Z)\big)=0$ for $i=1,2$,
\end{itemize}
the proposition follows.
\end{proof}

The following theorem shows that when $\h$ is  
solvable and $\beta(X)$ is a semisimple derivation of $\h$ for all $X\in\a$, 
then the total cohomology is preserved by elementary infinitesimal deformations.

\begin{theorem}\label{thm.total_deformation}
Let $\h$ be a solvable Lie algebra, $\a$ an abelian Lie algebra, 
$\a\oplus_{\alpha,\rho}\h$ an extension of $\a$ by $\h$ (necessarily solvable) and 
$\a\oplus_{\beta,0}\h_{\text{ab}}$ 
an elementary infinitesimal linear deformation of 
$\a\oplus_{\alpha,\rho}\h$ such that $\beta(X)$ is a
semisimple derivation of $\h$ for all $X\in\a$.
Then 
\[
\dim TH^p(\a\oplus_{\alpha+\beta,\rho}\h)=\dim TH^p(\a\oplus_{\alpha,\rho}\h).
\]
\end{theorem}

\begin{proof} Let  $\g_{\alpha,\rho}=\a\oplus_{\alpha,\rho}\h$
and $\g_{\alpha+\beta,\rho}=\a\oplus_{\alpha+\beta,\rho}\h$.
We proceed by induction on $\rank\,\beta=\dim(\Im\beta)$.
If $\rank\beta=1$, let $X_1\in\a$ such that $\beta(X_1)\ne0$
and let $\a_1=\K X_1$, $\a_2=\ker \beta$ and $\s=\a_2+\h\subseteq\g_{\alpha,\rho}$.
Clearly, $\s$ is an ideal of $\g_{\alpha,\rho}$
and if $D=\ad_{\g_{\alpha,\rho}}(X_1)|_{\s}\in\Der(\s)$, then we have 
\[
  \g_{\alpha,\rho}=D\ltimes\s.
\]
Let $S\in\End(\s)$ be defined by
\[
  S(A)=
\begin{cases}
 0, & \text{if }A\in\a_2; \\
 \beta(X_1)(A), & \text{if }A\in\h.
\end{cases}
\]
Since $\beta(X_1)\in\Der(\h)$ and it satisfies properties (3'') and (4''),
it follows that $S\in\Der(\s)$ and 
\[
 \g_{\alpha+\beta,\rho}=(D+S)\ltimes\s. 
\]
Properties (3'') and (4'') also imply 
that $S$ and $D$ commute with each other and hence
it follows from Corollary \ref{coro.Total_Ssimple} that 
$\dim TH^p(\g_{\alpha+\beta,\rho})=\dim TH^p(\g_{\alpha,\rho})$.

If $\rank\,\beta = n>1$, then 
choose a non trivial decomposition of $\a=\a_1+\a_2$ 
that yields a decomposition $\beta=\beta_1+\beta_2$ with
 $\rank\,\beta_i < n$ for $i=1,2$. 
It follows from Proposition \ref{prop.Elementary} and 
the induction hypothesis applied twice that
\[
\dim TH^p(\g_{\alpha+\beta,\rho})
 =\dim TH^p(\g_{\alpha+\beta_1,\rho}) 
  =\dim TH^p(\g_{\alpha,\rho}) 
\]
and this completes the proof.
\end{proof}

\section{The nilshadow of a solvable Lie algebra}\label{sec:nilshadow}

Let $\s$ be a solvable Lie algebra and let $\n$ be its nilradical.
If $X\in\s$, let $\text{ad}(X) = \ad(X)_s + \ad(X)_n$ be
the Jordan decomposition of $\ad(X)$ in $\text{End}(\s)$. 
It is well known that both $\text{ad}(X)_s$ and $\text{ad}(X)_n$ are derivations of $\s$
and therefore their image is contained in $\n$ (see for 
instance \cite[Lemma 4.2B]{H} and \cite[Corollary 2, Ch II.7]{J}).

Given a Cartan subalgebra $\h$ of $\s$ (a nilpotent self-normalizing subalgebra) it is not difficult 
to see that there exists a vector subspace $\a$ of $\h$ such that $\s = \a\oplus\n$ as vector spaces.
Moreover $\{\text{ad}(X)_s:X\in\a\}$ is an abelian subalgebra of $\text{Der}(\s)$
of semisimple operators vanishing on $\a$ and preserving $\mathfrak{n}$ (see \cite{DER}).
In fact, if we define $\beta:\a\to\text{Der}(\n)$ by
\[
 X \mapsto \ad(X)_s|_{\n},
\]
then $\beta$ is a Lie algebra homomorphism satisfying (3'') and (4'') and thus 
$\a\oplus_{\beta,0}\n_{\text{ab}}$ is an elementary infinitesimal linear deformation of $\s$.
On the other hand, the bracket
\begin{align*}
[X,Y]_{\nilsh}&=[X,Y]-\ad(X_\a)_s(Y)+\ad(Y_\a)_s(X), \\
           &=[X,Y]-\ad(X_\a)_s(Y_\n)+\ad(Y_\a)_s(X_\n),  
\end{align*}
where $X=X_\a+X_\n$ is the decomposition of $X$ with respect to $\s=\a+\n$,
defines on $\s$ a new Lie algebra structure
which is nilpotent and coincides on the nilradical of
$\s$ with the original one.
This new Lie algebra is the \emph{nilshadow} of $\s$ and we denote it by $\nilsh(\s)$.

The nilshadow was originally considered in the context of groups in \cite{AG}
and for Lie algebras can be found in \cite{B} and \cite{DER}.

\begin{remark}
The nilshadow constructed above depends on the choice of the Cartan subalgebra $\h$ and the subspace $\a$, 
but different choices yield isomorphic nilshadows.
\end{remark}

\begin{remark}
Since  $\dim\{\text{ad}(X)_s:X\in\a\}=\dim\a$, it is not difficult to see that 
the nilshadow described above is isomorphic to the nilradical of $\text{ad}(\a)_s\ltimes \s$.
\end{remark}

\begin{theorem}
 The nilshadow $\nilsh(\s)$ of a solvable Lie algebra 
$\s$ is an elementary linear deformation of $\s$ and in particular
\[
\dim TH^p(\s)=\dim H^p(\nilsh(\s)),
\]
for all $p=0,\dots,\dim\s$.
\end{theorem}

\begin{proof}
Let $\alpha$ and $\rho$ be such that $\s=\a\oplus_{\alpha,\rho}\n$.
 It is straightforward to see that the bracket $[X,Y]_{\nilsh}$ on 
$\s=\a\oplus\n$ is the same as in $\a\oplus_{\alpha-\beta,\rho}\n$.
As we already pointed out, 
$\a\oplus_{\beta,0}\n_{\text{ab}}$ is an elementary infinitesimal linear 
deformation of $\s$.
Now the theorem follows from Theorem \ref{thm.total_deformation}.
\end{proof}

We conclude with a remark about different linear subspaces of linear deformations 
associated to a solvable Lie algebra, all of which have the same total cohomology.

Let $\s$ be a solvable Lie algebra with nilradical $\n$ and let $r=\dim(\s/\n)$. 
Given an ideal $\h$ of $\s$ containing $\s'$ and 
a complementary subspace $\a$ of $\h$ (that we 
may assume is contained in a Cartan subalgebra of $\s$) let 
\[
 \b_{\h,\a}=\{D\in\Der(\h):[D,\ad(\a)]=0\text{ and }D|_{[\a,\a]}=0\}.
\]
It is clear that  $\b_{\h,\a}$ is a Lie subalgebra of $\Der(\h)$.
For any toral subalgebra $\mathfrak{t}$ of $\b_{\h,\a}$ and $\beta:\a\to\mathfrak{t}$
a linear map, it follows from \S\ref{sec:deformations} that 
\[
 \s_\beta=\a\oplus_{\alpha+\beta,\rho}\h
\]
is an elementary linear deformation of $\s$ and Theorem \ref{thm.total_deformation}
says that 
\[
\dim TH^p(\s)=\dim H^p(\s_\beta).
\]
Notice that one choice for $\mathfrak{t}$ is the subspace consisting  
of all the semisimple parts of $\ad(X)$ with $X\in\a$ (see \cite{DER}).
In this case, $\dim\mathfrak{t}\ge\min\{\dim\a,r\}$.
In other words we have shown that there exists
a linear space of Lie algebras 
\[
\{\s_\beta:\beta\in\Hom(\a,\mathfrak{t})\}
\]
with dimension at least $\dim\a\times\min\{\dim\a,r\}$
such that $\s_0=\s$ and  all its members have the same total cohomology.
In the particular case when $\h=\n$ this choice of $\mathfrak{t}$
leads to a linear space of Lie algebras of dimension at least $r^2$ 
containing the nilshadow of $\s$.


\end{document}